\documentclass[12pt]{amsart}
\usepackage{a4wide}

\usepackage[utf8]{inputenc}
\usepackage[english]{babel}
\usepackage{amsmath, amssymb, amsthm, bm, braket, array, mathrsfs, mathtools}
\usepackage{enumerate}
\usepackage{subfigure}
\usepackage[bookmarks]{hyperref}

\newtheorem{theorem}{Theorem}[section]
\newtheorem{proposition}[theorem]{Proposition}
\newtheorem{lemma}[theorem]{Lemma}

\theoremstyle{definition}

\newtheorem{problem}[theorem]{Problem}

\theoremstyle{definition}
\newtheorem{definition}[theorem]{Definition}

\newcommand{\N}{\mathbb{N}}

\newcommand{\Pfin}{\mathfrak P_{\mathrm{fin}}}
\newcommand{\Pboole}{\mathfrak P_{2}}
\newcommand{\Pthree}{\mathfrak P_{3}}
\newcommand{\Pn}{\mathfrak P_{n}}

\newcommand{\struct}[1]{\mathbb{#1}}
\newcommand{\clo}[1]{\mathcal{#1}}
\newcommand{\alg}[1]{\mathbf{#1}}

\newcommand{\algA}{\alg{A}}
\newcommand{\algB}{\alg{B}}

\newcommand{\algS}{\alg{S}}

\DeclareMathOperator{\pr}{pr}
\DeclareMathOperator{\CSP}{CSP}
\DeclareMathOperator{\PCSP}{PCSP}
\DeclareMathOperator{\Pol}{Pol}
\DeclareMathOperator{\Con}{Con}
\DeclareMathOperator{\Clo}{Clo}
\DeclareMathOperator{\Inv}{Inv}

\newcommand{\isaminorof}{\leq_{\mathrm{m}}}
\newcommand{\notisaminorof}{\nleq_{\mathrm{m}}}
\newcommand{\minorequivalentto}{\equiv_{\mathrm{m}}}

\DeclareMathOperator{\HSP}{HSP}
\DeclareMathOperator{\Ho}{H}
\DeclareMathOperator{\Su}{S}
\DeclareMathOperator{\Prd}{P}

\DeclareMathOperator{\HS}{HS}

\newcommand{\tp}[1]{\mathbf{#1}}

\newcommand{\minority}{\mathtt{min}}

\begin{document}
\title{Mal'cev clones over a three-element set up to minor-equivalence}

\author[S.~Fioravanti]{Stefano Fioravanti}
\address{Department of Algebra\\ Charles University\\Prague\\ Czechia
}
\email{stefano.fioravanti66@gmail.com}

\author[M.~Kompatscher]{Michael Kompatscher}
\address{Department of Algebra\\ Charles University\\Prague\\ Czechia}
\email{kompatscher@karlin.mff.cuni.cz}

\author[B.~Rossi]{Bernardo Rossi}
\address{Siena, Italy}
\email{bernardo.rossi96@gmail.com} 

\author[A.~Vucaj]{Albert Vucaj}
\address{Poggibonsi, Italy}
\email{albo.bradbury@gmail.com}

\thanks{This paper was supported by Charles University under grant number PRIMUS/24/SCI/008. Michael Kompatscher is further funded by the Czech Science Foundation
(GA\v{C}R grant no. 25-16324S), and Charles University Research Center program No. UNCE/24/SCI/022.}

\begin{abstract}
    We classify all Mal'cev clones over a three-element set up to minion homomorphisms. This is another step toward the complete classification of three-element relational structures up to pp-constructability.
    We furthermore provide an alternative proof of Bulatov's result that all Mal'cev clones over a three-element set have an at most 4-ary relational basis.
\end{abstract}



\maketitle

\section{Introduction}\label{sec:Intro}
This article is dedicated to the study of clones over finite sets, a central topic in universal algebra that traces its origins back to Post’s classification of clones on a two-element set~\cite{Post}. It is well known that clones on a finite set with $n$ elements form a lattice when ordered by inclusion. In the following, we will denote this lattice by $\mathfrak{L}_n$. 

By Post's result, we know that $\mathfrak{L}_2$ is countable; however, Janov and Muchnik proved that on every finite set with at least three elements there are continuum many clones~\cite{3elem}. This finding suggests that a classification of $\mathfrak{L}_n$ \`{a} la Post becomes significantly more challenging for $n\geq 3$. In fact, up to date, no complete description of $\mathfrak{L}_3$ is known. Nevertheless, some remarkable progress has been made on the description of several significant parts of $\mathfrak{L}_3$, such as the classification of all 18 maximal~\cite{JabMaximal} and 84 minimal clones~\cite{Csakany84}. Notably, these results were obtained prior to Rosenberg’s groundbreaking classification results on maximal respectively minimal clones on general finite sets~\cite{RosenbergMaximal} and~\cite{Rosenberg}. One of the challenges in achieving a full description of $\mathfrak{L}_3$ lies in the fact that, with the exception of the clone of linear operations, every maximal clone generates an uncountable ideal~\cite{DemetrHannak, Marchenkov}. In 2015, Zhuk provided a classification of all clones of self-dual operations on a three-element set~\cite{Zhuk15}, which, to date, remains the only complete description of such an uncountable maximal ideal.

In this paper, we are interested in the subset of $\mathfrak{L}_3$ consisting of all \emph{Mal'cev clones}. A clone is said to be a Mal'cev clone if it contains a \emph{Mal'cev operation}, i.e., a ternary operation $m$ satisfying the identity $y \approx m(y,x,x) \approx m(x,x,y)$. In~\cite{BulatovImportante} Bulatov provided a complete classification of all 1129 Mal'cev clones in $\mathfrak L_3$ and described each of them by an at most 4-ary relational basis. In~\cite{KS-parallelogramterms} Kearnes and Szendrei observed that, for every $k \in \N$, there are also only finitely many clones in $\mathfrak{L}_3$ with a $k$-edge term (although without explicitly classifying them). Our first contribution is a simplified proof of Bulatov's result that uses the methods introduced in~\cite{KS-parallelogramterms}.

However, for many purposes, it is convenient to study clones up to a coarser order than inclusion. In particular, many properties of a clone are determined by the identities it satisfies. This is captured by the (pre-)ordering of clones by the existence of a \emph{clone homomorphism}, that is, a map between clones that preserves identities. This order is of particular interest in classic universal algebra, as, for two algebraic structures $\algA, \algB$ (without constant symbols), there is a clone homomorphism from the clone of term operations $\Clo(\algA)$ to the clone of term operations $\Clo(\algB)$ if and only if the variety generated by $\algA$ has an interpretation in the variety generated by $\algB$. Thus, studying clones up to clone homomorphism essentially corresponds to classifications of varieties in the interpretability lattice introduced by Garcia and Taylor~\cite{garciaTaylor}.

Motivated by the study of fixed-template constraint satisfaction problems (CSPs), an even coarser order was introduced by Barto, Opr\v{s}al, and Pinsker via the notion of \emph{pp-constructability}~\cite{wonderland}. Pp-constructions form a unifying framework that describes several different types of reductions between CSPs. More precisely, if a finite structure $\struct{A}$ pp-constructs a finite structure $\struct{B}$, then there exists a log-space reduction from $\CSP(\struct{B})$ to $\CSP(\struct{A})$. 

The algebraic counterpart of pp-constructability is a weakened version of clone homomorphism, referred to in the literature as \emph{minion homomorphism} or \emph{minor-preserving map}. For any pair of finite relational structures $\struct{A}$ and $\struct{B}$, $\struct{A}$ pp-constructs $\struct{B}$ if and only if there exists a minion homomorphism from 
$\Pol(\struct{A})$ to $\Pol(\struct{B})$~\cite{wonderland}. Here $\Pol(\struct{A})$ denotes the \emph{polymorphism clone} of $\struct{A}$, that is, the clone consisting of all operations that preserve the relations of $\struct{A}$. Thus, the polymorphism clone, up to minion homomorphism, completely determines the complexity of $\CSP(\struct{A})$. This is a fundamental fact used in the ``universal-algebraic approach'' to CSPs, which leads to numerous results, among which the most famous is probably the P/NP-c dichotomy for all finite CSPs, independently proved by Bulatov~\cite{BulatovFVConjecture} and Zhuk~\cite{ZhukFVConjecture, ZhukDichotomy}. More recently, pp-constructions were also observed to describe reductions between \emph{promise constraint satisfaction problems, $\PCSP(\struct{A},\struct{B})$}; also in this setting, pp-constructions between finite PCSP-templates correspond to minion homomorphisms between their \emph{polymorphism minons} $\Pol(\struct{A},\struct{B})$, hence the name~\cite{jakubPCSP}.

The main difference between clone and minion homomorphism lies in the fact that clone homomorphisms preserve \emph{all} identities, while minion homomorphisms only preserve so-called \emph{minor identities} or \emph{height 1 identities}, i.e., identities that can be stated without using composition and projections. We will define this notion precisely in Definition~\ref{def:minor}.

If we preorder the class of all clones over some finite set by the existence of minion homomorphisms and factor by the induced equivalence relation $\minorequivalentto$, we obtain the poset $\Pfin$. In particular, we say that two clones $\clo{C}$ and $\clo{D}$ are \emph{minor-equivalent} if $\clo{C} \minorequivalentto \clo{D}$. If we restrict ourselves to clones over $n$-element sets, we obtain $\Pn = \mathfrak L_n/_{\minorequivalentto}$. A systematic study of $\Pfin$ and $\Pn$ was initiated in the Ph.D. theses of Starke~\cite{FlorianThesis} and the fourth author~\cite{AlbertThesis}. Although we do not know whether $\Pfin$ is a lattice, we know that it is a $\wedge$-semilattice without atoms and a unique coatom~$\overline{\clo{I}_2}$. Moreover, $\Pfin$ contains infinite descending chains and infinite antichains (see, e.g.,~\cite{cyclesBodirskyStrakeVucaj}). A full description of all submaximal elements of $\Pfin$ was recently provided in~\cite{meyerSubmaximal}, generalizing two previous results that were limited to $\Pthree$~\cite{VucajZhuk} and polymorphism clones of finite digraphs~\cite{BodirskyStarke}. A complete description of $\Pboole$ was provided in~\cite{albert}. This was also used as a basis for the classification of multisorted Boolean clones determined by binary relations~\cite{BartoKapytka} and products of Boolean clones~\cite{RadekOlsak} up to minion homomorphism.

The complete structure of the poset $\Pthree$ remains yet unknown; only partial results have been proven so far. By~\cite{VucajZhuk}, there are exactly three submaximal elements, which are the equivalence classes of the clones $\clo{B}_2 = \Pol(\mathbb B_2)$, $\clo{C}_2 = \Pol(\mathbb C_2)$, and $\clo{C}_3 = \Pol(\mathbb C_3)$. Here, $\mathbb C_n$ denotes the cyclic digraph of $n$-vertices with constants, and $\mathbb B_2 = (\{0,1\}; \{0,1\}^2\setminus \{(0,0)\}, \{0\},\{1\})$. 

Furthermore, the ideal of clones of self-dual operations in $\mathfrak L_3$, i.e., the ideal generated by $\clo{C}_3$ in $\mathfrak L_3$, was completely described up to minion homomorphisms in~\cite{BodirskyVucajZhuk}.

The main result of our paper, Theorem~\ref{teor:each_malcev_clone_in_one_class}, is a classification of all Mal'cev clones on a three-element set up to minor-equivalence. In fact, our classification applies more generally to all clones with a \emph{quasi Mal'cev operations} $m(x,x,y)\approx m(y,x,x)\approx m(y,y,y)$. We show that there are exactly 10 equivalence classes; their order is depicted in Figure~\ref{figure:main_figure}.

Our result can be seen as a significant step towards a full description of $\Pthree$. Indeed, as was observed in \cite[Theorem 2.13]{VucajZhuk}, every element in $\Pfin$ that does not have a quasi-Malcev operation lies below $\mathcal B_2$. In particular, this is also true in $\Pthree$. Therefore, by our results, it would suffice to describe the downset of $\mathcal B_2$ in order to complete the classification of $\Pthree$.

\subsection*{Organization of the paper} In Section~\ref{sec:Preliminaries}, we introduce some necessary definitions and notations. Section~\ref{sec:malcev} recalls some results from~\cite{KS-parallelogramterms}, which we then use in Section~\ref{sec:reprovingBulatov} to provide an alternative proof of the main result in~\cite{BulatovImportante}. In the remaining parts of the paper, we focus on the classification of Mal’cev clones on a three-element set, up to minion homomorphism. We distinguish several subcases:

\begin{enumerate}[(i)]
    \item In Section~\ref{sec:selfdual} we identify classes that are covered by the description of $\Pboole$ in \cite{albert}, and the classification of clones of self-dual operations in~\cite{BodirskyVucajZhuk}. \label{itm:Intro1}
    \item In Section~\ref{sec:SDmeet} we classify idempotent clones with a majority operation. \label{itm:Intro2} \label{itm:Intro3}
    \item All remaining cases are discussed in Section~\ref{sec:NotSDmeet}.\label{itm:Intro4} 
\end{enumerate}
We then present the main result of this paper in Section~\ref{sec:Main}, including a comprehensive overview of the 10 classes in Figure~\ref{figure:main_figure} and Figure~\ref{fig:Separations}. In Section~\ref{sec:future} we discuss some directions for future research. The appendix contains a link to a database that, among other information, describes the $\minorequivalentto$ equivalence class for each of the 1129 clones in Bulatov's classification.

\section{Preliminaries}\label{sec:Preliminaries}
In this section, we introduce some necessary notation, key definitions, and foundational results from the literature. For further background on clones we refer to \cite{PosKal79,Bod21}. For a more comprehensive presentation of the notions and results in universal algebra that we use, we refer to \cite{McKMcnTay88} and \cite{BS}.

\subsection*{Basic notation} We use boldface letters for tuples and index their entries by subscripts, e.g., $\boldsymbol{a}=(a_1,\dots,a_n)$. A tuple is \emph{injective} if all its components are distinct, i.e., if $|\{a_1,\dots,a_n\}|=n$.  We write $[n]$ for the set $\{1,\dots,n\}$. We are going to denote algebraic structures by capital boldface characters (e.g., $\alg{A} = (A;f_1,f_2,f_3)$) and relational structures by capital `blackboard' characters (e.g., $\struct{A} = (A;R_1,R_2)$).

We note that blackboard characters used for relational structures may overlap with the standard notation for the field of integers modulo a prime $p$, denoted by $\mathbb Z_p$. To avoid ambiguity, we will explicitly clarify when the latter interpretation is intended.

\subsection{Operation clones}\label{sec:operationClones}
For a set $A$ and $n \in {\mathbb N}$, we define ${\mathcal O}^{(n)}_A = \{f\mid f\colon A^n\to A\}$ and ${\mathcal O}_A = \bigcup_{n \in {\mathbb N}} {\mathcal O}^{(n)}_A$. For each $n\geq 1$ and $i\in\{1,\dots,n\}$, the $n$-ary \emph{projection} to the $i$-th coordinate is the
operation $\pr^{n}_{i}$ defined for all $a_1,\dots,a_n\in A$ by $(a_1,\dots,a_n)\mapsto a_i$. 
The \emph{composition of $f \colon A^n \to A$ and $g_1,\dots, g_n \colon A^m \to A$} is defined as the $m$-ary operation
\[
(f \circ (g_1,\dots, g_n))(a_1,\dots, a_m) = f(g_1(a_1,\dots, a_m),\dots, g_n(a_1,\dots, a_m)),\]
\noindent for $a_1,\dots, a_m\in A$.

A \emph{clone over $A$} is a subset $\clo{C}\subseteq\mathcal O_A$ which is closed under composition of operations and contains all projections. If $F \subseteq {\mathcal O}_A$, then $\Clo(F)$ denotes the \emph{clone generated by $F$}, i.e., the smallest clone that contains~$F$. For an algebraic structure $\algA$, we also write $\Clo(\algA)$ for the clone generated by the basic operations of $\algA$. This is commonly referred to as the \emph{clone of term operations} of $\algA$.

A clone $\mathcal C$ is called \emph{idempotent} if all its operations are idempotent, i.e., $f(x,x,\ldots,x) \approx x$ holds for all $f\in \mathcal C$. The \emph{idempotent reduct} of a clone is the subclone of $\mathcal C$ that consists of all idempotent operations in $\mathcal C$.
We say that a clone is a \emph{Mal'cev clone} if it contains a \emph{Mal'cev operation}, i.e. a ternary operation $d$ satisfying the identities $d(y,x,x) \approx d(x,x,y) \approx y$.

\subsection{Algebraic structures}

Every clone $\mathcal C$ on a set $A$ is equal to the
clone of term operations of some algebraic structure $\algA$ (for example, the algebra $\algA_\mathcal C$ whose
basic operations are exactly the elements of $\mathcal C$). As some parts of this paper are more readily explained in terms of algebraic structures rather than clones, we will switch between these two languages repeatedly. 

We assume that the reader is familiar with standard notions in universal algebra. In particular, we write $\alg{S} \leq \algA$, if $\algS$ is a \emph{subalgebra} of $\algA$. We denote the lattice  of \emph{congruences} of $\algA$ by $\Con(\algA)$, and write $\algA/\theta$ for the \emph{quotient} of $\algA$ by a congruence $\theta \in \Con(\algA)$. We use the notation $0_\algA = \{(x,x) \mid x \in A\}$ for the smallest congruence and $1_\algA = A\times A$ for the biggest congruence of $\algA$. An algebra $\algA$ is called \emph{subdirectly irreducible}, if it has a smallest non-trivial congruence $\mu \neq 0_\algA$. In this case $\mu$ is called the \emph{monolith} of $\algA$. We write $\prod_{i\in I}\algA_i$ for the \emph{product} of a family of algebras $(\algA_i)_{i\in I}$ (of the same language). For products of finitely many algebras $\algA_1, \ldots, \algA_n$, we also use the notation $\algA_1 \times \cdots \times \algA_n$.

For a given class of algebras $\mathcal K$ (of the same language), we denote its closure under products by $\Prd(\mathcal K)$, its closure under homomorphic images by $\Ho(\mathcal K)$, and its closure under subalgebras by $\Su(\mathcal K)$. Note that $\Ho(\mathcal K)$ is equal to the closure of $\mathcal K$ under quotients and isomorphisms. Any class of algebras that is closed under $\Prd$, $\Ho$ and $\Su$ is called a \emph{variety}. It is well-known that this is equivalent to $\mathcal K = \HSP(\mathcal K)$~\cite[Theorem 9.5]{BS}. By Birkhoff's famous HSP theorem~\cite{Bir-On-the-structure}, this is further equivalent to $\mathcal K$ being the class of models of an equational theory.

Let us say an algebra $\algA$ is \emph{idempotent} if $\Clo(\algA)$ is idempotent. Since we will mostly deal with idempotent clones, we would like to point out the useful observation that in idempotent algebras, every class of every congruence relation is a subalgebra of $\algA$.

For more background on universal algebra, we refer to~\cite{BS}.

\subsection{Relational structures and cores}
A $\tau$-structure is a relational structure over some signature $\tau$. Let $\struct{A}$ and $\struct{B}$ be two relational $\tau$-structures. A map $h\colon A\to B$ is a \emph{homomorphism} if, for every $R\in\tau$ and each $(a_1,\dots, a_n)\in A^n$, we have
\begin{equation*}
    (a_1,\dots,a_n)\in R^\struct{A} \Rightarrow (h(a_1),\dots,h(a_n))\in R^\struct{B}.
\end{equation*}
We also use the notation $h\colon \struct{A} \to \struct{B}$ if $h$ is a homomorphism from $\struct{A}$ to $\struct{B}$. 
We say that $\struct{A}$ and $\struct{B}$ are \emph{homomorphically equivalent} if there is a homomorphism from $\struct{A}$ to $\struct{B}$ and one from $\struct{B}$ to $\struct{A}$. 
An \emph{endomorphism of} $\struct{A}$ is a homomorphism from $\struct{A}$ to itself. An \emph{isomorphism} between $\struct{A}$ and $\struct{B}$ is a bijective homomorphism $h$ such that its inverse is also a homomorphism. An \emph{automorphism of} $\struct{A}$ is an isomorphism from $\struct{A}$ onto itself. A finite structure $\struct{A}$ is a \emph{core} if every endomorphism of $\struct{A}$ is also an automorphism. We say that $\struct{C}$ is a \emph{core of} $\struct{A}$ if $\struct{C}$ is a core and is homomorphically equivalent to $\struct{A}$. It is well established that every finite relational structure has a core that is unique up to isomorphism, thus, it makes sense to refer to \emph{the} core of a relational structure.

\subsection{The Inv-Pol Galois connection and pp-formulas}\label{sec:InvPol}
We say that an operation $f\colon A^n \to A$ \emph{preserves a relation} $R\subseteq A^k$ if, for any $\boldsymbol{a}_1,\dots, \boldsymbol{a}_n \in R$, it follows that $f(\boldsymbol{a}_1,\dots,\boldsymbol{a}_n) \in R$. Here, $f(\boldsymbol{a}_1,\dots,\boldsymbol{a}_n) \in A^k$ is computed by coordinate-wise evaluation. Equivalently, we also say that $R$ is \emph{invariant under} $f$. If $f$ preserves every relation from a set $\Gamma$, then $f$ is called a \emph{polymorphism of} $\Gamma$. We denote the set of all polymorphisms of $\Gamma$ by $\Pol(\Gamma)$. As $\Pol(\Gamma)$ forms a clone (cf. \cite[Chapter E.]{PosKal79}), it is commonly referred to as the \emph{polymorphism clone of} $\Gamma$. When working with relational structures, we also use the notation $\Pol(\struct{A})$ for $\Pol(\Gamma)$, where $\struct{A}=(A;\Gamma)$. 
    
Note that $\Pol(\struct{A})$ is the set consisting of all homomorphisms $f\colon \struct{A}^n\to \struct{A}$, for each $n\geq 1$. Here $\struct{A}^n$ denotes the $n$-th cartesian power of $\struct{A}$, i.e. the structure of domain $A^n$ with the same signature $\tau$ as $\struct{A}$ such that for all $R \in \tau$ of arity $k$:
\[
R^{\struct{A}^n} = \{(\boldsymbol{a}^1,\dots,\boldsymbol{a}^k) \in (A^n)^{k} \mid R^{\struct{A}}(a_j^1,\dots,a_{j}^k) , 1 \leq j \leq n\}.
\]
For any set $F$ of operations over $A$, we denote by $\Inv(F)$ the set of all relations that are invariant under every $f\in F$; and for an algebraic structure $\algA = (A;F)$, we write $\Inv(\algA) = \Inv(F)$. Note that $\Inv(F)$ consists of all subuniverses $R \leq \algA^n$, for each $n\geq 1$. We say that $\Gamma$ is a \emph{relational basis} of an algebra $\algA$ if $\Clo(\algA) = \Pol(\Gamma)$.

If $A$ is a \emph{finite} set, then $\Clo(\algA) = \Pol(\Inv(\algA))$ holds; hence, every clone on $A$ is the polymorphism clone of some relational structure~\cite{BoKaKoRo,Geiger}. There also exists a characterization of $\Inv(\mathbb A)$ in terms of pp-definitions. 
A \emph{primitive positive (pp-) formula over $\tau$} is a first-order formula which only uses relation symbols in $\tau$, equality, conjunction, and existential quantification. If $\struct{A}$ is a $\tau$-structure over a set $A$ and $\phi(x_1,\dots,x_n)$ is a $\tau$-formula with free-variables $x_1,\dots,x_n$, then 
$\{(a_1,\dots,a_n) \in A^n\mid \struct{A} \models \phi(a_1,\dots,a_n)\}$
is called \emph{the relation defined by $\phi$ in $\struct{A}$}. 
Any relation defined by some primitive positive formula over $\struct{A}$ is said to be \emph{pp-definable} in $\struct{A}$. If $\struct{A}$ and $\mathbb B$ are relational structures on the same domain, we say that $\struct{A}$ \emph{pp-defines} $\struct{B}$ if every relation in $\struct{B}$ is pp-definable in $\struct{A}$. It is easy to check that $\Inv(\Pol(\struct{A}))$ contains the set of all relations that are pp-definable from $\struct{A}$. The following result shows that if $\struct{A}$ has a finite domain, then 
$\Inv(\Pol(\struct{A}))$ is equal to the set of pp-definable relations:

\begin{theorem}[\cite{BoKaKoRo,Geiger}]\label{thm:DefAndCloneInclusion}
Let $\struct{A}$ and $\struct{B}$ be structures on the same finite set $A$. Then $\struct{A}$ pp-defines $\struct{B}$ if and only if $\Pol(\struct{A})\subseteq\Pol(\struct{B})$.
\end{theorem}

Note that Theorem \ref{thm:DefAndCloneInclusion} implies that $\Pol(\struct{A})$ is idempotent if and only if $\struct{A}$ pp-defines all unary singleton relations. In this case $\struct{A}$ is also a core, as the endomorphisms of $\struct{A}$ only consist of the identity on $A$.

\subsection{Minors and minion homomorphisms} \label{sec:ourposet}

In this subsection we introduce the notions of minor and minion homomorphism.
\begin{definition}\label{def:minor}
Let $f\colon A^n \to A$ be an $n$-ary operation and let $\sigma \colon [n] \to [r]$. We denote by $f_\sigma$ the following $r$-ary operation 
\begin{equation*}
f_\sigma(x_1,\dots,x_{r}) = f(x_{\sigma(1)},\dots,x_{\sigma(n)}).
\end{equation*}
Any operation of the form $f_\sigma$, for some $\sigma \colon [n]\to [r]$, is called a \emph{minor} of $f$.
\end{definition}

In other words, the minors of an operation $f$ are exactly those functions that can be obtained by permuting/identifying variables and adding dummy variables.

\begin{definition}\label{def:minorpreservingmaps}
Let $\clo{A}$ and $\clo{B}$ be clones over some finite (possibly distinct) sets, and let $\xi \colon \clo{A} \to \clo{B}$ be a mapping that preserves arities. We say that $\xi$ is a \emph{minion homomorphism} if
\begin{equation*}
\xi(f_\sigma) = \xi(f)_\sigma
\end{equation*}
for any $n$-ary operation $f \in \clo{A}$ and $\sigma\colon [n] \to [r]$. 

We write $\clo{C}\isaminorof\clo{D}$ if there exists a minion homomorphism
from $\clo{C}$
to $\clo{D}$. 
If $\clo{C}\isaminorof\clo{D}$ and
$\clo{D}\isaminorof\clo{C}$,
we say that $\clo{C}$ and $\clo{D}$ are
\emph{minor-equivalent} and write
$\clo{C}\minorequivalentto\clo{D}$. For any clone $\clo{C}$, we also write $\overline{\clo{C}}$ to denote the $\minorequivalentto$-class of $\clo{C}$.
\end{definition}

Clearly $\clo{C}\subseteq\clo{D}$ implies $\clo{C}\isaminorof\clo{D}$. Let us remark that minion homomorphisms can also be defined for \emph{minions} as defined in~\cite{jakubPCSP}, and that the image of a clone under a minion homomorphism does not need to be a clone.

A \emph{minor identity} or \emph{height 1 identity} is a formal expression of the form 
\[f(x_{\sigma(1)},\dots,x_{\sigma(n)}) \approx g(x_{\pi(1)},\dots,x_{\pi(m)}),\] 
for some function symbols $f$, $g$, and some $\sigma\colon [n] \to [r]$ and $\pi\colon [m] \to [r]$. A \emph{minor condition} $\Sigma$ is a finite set of minor identities. 
We say that a set of operations $F$ satisfies a minor condition $\Sigma$ (and write $F\models \Sigma$) if we can assign each function symbol in $\Sigma$ to an operation in $F$ in such a way that all the (universally quantified) identities in $\Sigma$ are satisfied. For instance, the minor condition $\Sigma = \{f(x,y) \approx f(y,x) \}$ is satisfied by $F$ if $F$ contains a binary commutative operation. 
Every minion homomorphism preserves minor conditions. By a compactness argument outlined in \cite{wonderland}, there is a minion homomorphism from $\clo{C}$ to $\clo{D}$ if and only if $\clo{C} \models \Sigma \Rightarrow \clo{D} \models \Sigma$ holds for every minor condition $\Sigma$. 

We refer to the partially ordered set of $\minorequivalentto$-classes of clones over a finite set as the \emph{minion homomorphism poset} and denote it by $\Pfin$. We also consider the subposet of $\Pfin$ that consists of $\minorequivalentto$-classes containing a representative on an $n$ element set, which we denote by $\Pn$. More precisely,
\begin{align*}
    \Pfin &=(\{\overline{\clo{C}}\mid \clo{C}\text{ is a clone over some finite set}\}; \isaminorof), \text{ and}\\
    \Pn &=(\{\overline{\clo{C}}\mid \clo{C}\text{ is a clone over } \{0,\dots,n-1\}\}; \isaminorof),
\end{align*}
where $\overline{\clo{C}}\isaminorof\overline{\clo{D}}$ if and only if $\clo{C}\isaminorof\clo{D}$. 

We remark that a clone $\clo{C}$ on a set $A$ is minor equivalent to the clone $\clo{C^*}$ on the set $A \cup \{b\}$ that is generated by all functions $f^* \colon (A \cup \{b\})^n \to (A \cup \{b\})$ with $f^*|_{A^n} = f \in \clo{C}$, and $f^*(x_1,\ldots,x_n) = b$ if $x_i = b$ for some $i$ (both the map $f \mapsto f^*$, and $g \to g|_A$ can easily be seen to be minion homomorphism from $\clo{C}$ to $\clo{C^*}$ and vice-versa). This implies that, for every $n$, there is an order-preserving embedding of $\Pn$ into $\mathfrak P_{n+1}$.

Note that the top element of $\Pfin$ is the equivalence class of any clone that contains a constant operation; the bottom element of $\Pfin$ is the equivalence class containing the clone of projections on any set with at least two elements. These two classes are also the top and bottom elements of $\Pthree$.

\subsection{Pp-constructions and important minor conditions}
Let $\struct{A}$ and $\struct{B}$ be finite relational structures. We say that $\struct{B}$ is a \emph{pp-power} of $\struct{A}$ if $\struct{B}$ is isomorphic to a structure $\struct{P}=(P;\Gamma)$ such that
\begin{enumerate}
    \item $P=A^n$, for some $n\geq 1$,
    \item $\Gamma$ is a list of relations on $P$ which are pp-definable from $\struct{A}$.
\end{enumerate}
To be more precise, in (2), we require that any $k$-ary relation $R\subseteq (A^n)^k$ in $\Gamma$ is pp-definable from $\struct{A}$ when regarded as a $kn$-ary relation on $A$.

\begin{definition}
Let $\struct{A}$ and $\struct{B}$ be finite relational structures. We say that $\struct{A}$ \emph{pp-constructs} $\struct{B}$ if $\struct{B}$ is homomorphically equivalent to a pp-power of $\struct{A}$.
\end{definition}
In analogy to Theorem~\ref{thm:DefAndCloneInclusion}, it was shown in~\cite{wonderland} that pp-constructability between finite relational structures translates into the notion of \emph{minion homomorphism} between the corresponding polymorphism clones:

\begin{theorem}[\cite{wonderland}]\label{thm:minhomAndPPConAreTheSameThing} Let $\struct{A}$ and $\struct{B}$ be finite relational structures. The following are equivalent:
\begin{enumerate}
     \item $\struct{A}$ pp-constructs $\struct{B}$,
    \item $\Pol(\struct{A})\isaminorof\Pol(\struct{B})$,
    \item If $\Pol(\struct{A})$ satisfies a minor condition $\Sigma$, then $\Pol(\struct{B})$ satisfies $\Sigma$.
\end{enumerate}
\end{theorem}

Thus, the order $\Pfin$ is anti-isomorphic to finite relational structures ordered by pp-constructability (in fact, this is how $\Pfin$ is defined in~\cite{albert}).

The next proposition shows that, in the context of classifying clones on finite sets up to minor-equivalence, we can restrict our focus to idempotent clones:

\begin{proposition}[\cite{wonderland}]\label{prop:idempotent}
For a finite relational structure $\mathbb A$, let $\mathbb A'$ denote its core, extended by all singleton unary relations. Then $\mathbb A$ pp-constructs $\mathbb A'$ and vice-versa. As a consequence, $\clo{C} = \Pol(\struct{A})$ is minor-equivalent to the idempotent clone $\clo{C}' = \Pol(\struct{A})'$.
\end{proposition}

By Theorem~\ref{thm:minhomAndPPConAreTheSameThing} (3), whenever $\struct{B}$ is not pp-constructable from $\struct{A}$, this is witnessed by a minor condition $\Sigma$. In the following definition, we introduce some minor conditions that are well known in the literature and that, together with those in Definition~\ref{def:MinorConditions2}, will be used to prove inequalities of the form $\clo{C}\notisaminorof\clo{D}$ in Section~\ref{sec:Main} (see Figure~\ref{fig:Separations}).

\begin{definition}\label{def:MinorConditions} 
We define the following minor conditions:
\begin{itemize}
\item A $n$-ary operation $c$ is called a \emph{$n$-cyclic operation} if it satisfies
\[
c(x_1,x_2,\dots,x_n)\approx c(x_2,\dots,x_n,x_1).
\]
\item A $n$-ary operation $f$ is a \emph{symmetric operation} if it satisfies
\[
    f(x_1, x_2,\dots, x_n) \approx f(x_{\pi(1)}, x_{\pi(2)},\dots,x_{\pi(n)}),
\]
for every permutation $\pi$ of $\{1,2,\dots,n\}$.
\item A ternary operation $m$ is a \emph{quasi minority operation} if it satisfies
\[
m(x,y,y)\approx m(y,x,y)\approx m(y,y,x)\approx m(x,x,x).
\]
\item A ternary operation $m$ is a \emph{quasi majority operation} if it satisfies
\[
m(x,y,y)\approx m(y,x,y)\approx m(y,y,x)\approx m(y,y,y).
\]
\end{itemize}
As usual in the literature, a \emph{minority (majority) operation} is an idempotent quasi minority (majority) operation.
\end{definition}

\section{Mal'cev algebras and commutator theory} \label{sec:malcev}

We recall some well-known facts about Mal'cev algebras and their congruences.

\begin{definition}A binary relation $R \subseteq A \times B$ has the \emph{parallelogram property} if $R(a,c),R(a,d),R(b,c)$ implies $R(b,d)$. An $n$-ary relation $R \subseteq \prod_{i=1}^n A_i$ has the \emph{parallelogram property} if, for each subset $\emptyset \neq I \subset [n]$, $R$ has the parallelogram property considered as a binary relation between $\prod_{i \in I}A_i$ and $ \prod_{j \notin I}A_j$.
\end{definition}

Whenever $R \leq \prod_{i=1}^n \alg{A}_i$ for some algebras $\alg{A}_i$ from a variety with a Mal'cev term, then $R$ has the parallelogram property. In fact, varieties with a Mal'cev term can be characterized by this property (see, e.g., \cite[Theorem 3.7.]{KS-parallelogramterms}). In particular, relations $R \leq \alg{A}^n$ for a Mal'cev algebra $\alg{A}$ always have the parallelogram property.

Next, we introduce some notions from commutator theory:

\begin{definition} \label{def:centralize}
Let $\algA$ be an algebra with Mal'cev term $d(x,y,z)$, and let $\alpha, \beta \in \Con(\algA)$. We say:
\begin{itemize}
\item \emph{$\alpha$ centralizes $\beta$} if $\algA$ preserves the 4-ary relation
$$T^\algA(\alpha,\beta) = \{(x,y,z,d(x,y,z)) \mid x \, \alpha \, y \, \beta \, z)\},$$
\item $\alpha$ is \emph{Abelian} if $\alpha$ centralizes itself, i.e. $\algA$ preserves
$$T^\algA(\alpha) = \{(x,y,z,d(x,y,z)) \mid x \, \alpha \, y \, \alpha \, z)\},$$
\item $\algA$ is \emph{Abelian} if the full congruence $1_A$ is Abelian, i.e. $\algA$ preserves
$$T^\algA = \{(x,y,z,d(x,y,z)) \in A^4 \mid x,y,z \in A\},$$
\item the \emph{centralizer of $\alpha$}, denoted by $(0_A:\alpha)$, is the largest $\delta \in \Con(\algA)$ such that $\alpha$ centralizes~$\delta$. 
\end{itemize}
\end{definition}
These notions are usually defined in a different way in the literature, which is based on the so-called `term-condition' \cite[Chapter 4]{FreeseMcKenzie}. However, in the special case of Mal'cev algebras, they all agree with Definition~\ref{def:centralize}.

\begin{lemma} \label{lemma:centerwelldef}
Let $\algA$ be an algebra with Mal'cev term $d$ and $\alpha,\beta \in \Con(\algA)$. Then:
\begin{enumerate}
\item All notions introduced in Definition~\ref{def:centralize} are well-defined and equivalent to the corresponding definitions in~\cite{FreeseMcKenzie}.
\item If $\alpha$ centralizes $\beta$, then $T^\algA(\alpha,\beta)$ does not depend on the choice of the Mal'cev term $d \in \Clo(\algA)$.
\item Assume that $\alpha\leq \beta$ and that $\alpha$ centralizes $\beta$, and let $a\in A$. Then $([a]_\alpha,d) = ([a]_\alpha,x-y+z)$, for some Abelian group $+$ on $[a]_\alpha$, and all algebras $([b]_\alpha,d)$ with $b \, \beta \, a$ are isomorphic to $([a]_\alpha,d)$.
\end{enumerate}
\end{lemma}

\begin{proof}
It follows from \cite[Exercise 6.7]{FreeseMcKenzie} that our definition of $\alpha$ centralizes $\beta$ is equivalent to the standard definition based on the term condition; for a proof, see, e.g., \cite[Lemma 2.4.]{AichingerMayrpqextenstions}. As a consequence, all other notions in Definition~\ref{def:centralize} are also equivalent to the standard definitions.
For (2), note that $(x,y,y,x), (y,y,y,y), (y,y,z,z) \in T^\algA(\alpha,\beta)$ holds for all $x \, \alpha \, y \, \beta \, z$. But then their image under any Mal'cev term $d'$ must also be in $T^\algA(\alpha,\beta)$, i.e., $(x,y,z,d'(x,y,z)) \in T^\algA(\alpha,\beta)$. Thus, $d'$ agrees with $d$ on such tuples. For point (3), note that the congruence class $[a]_\alpha$ is invariant under $d$, since $d$ is idempotent. If we define $+$ by $x+y = d(x,a,y)$, then it is a straightforward exercise to prove that $+$ is an Abelian group with neutral element $a$, and $([a]_\alpha,d) = ([a]_\alpha,x-y+z)$ (cf. \cite[Lemma 5.6]{FreeseMcKenzie}). Similarly, the preservation of $T^\algA(\alpha,\beta)$ can be used to prove that the map $x \to d(x,a,b)$ is an isomorphism from $([a]_\alpha,d)$ to $([b]_\alpha,d)$ (see \cite[Exercise 5.5]{FreeseMcKenzie}).
\end{proof} 

We next define a few useful notions for relations with the parallelogram property:

\begin{definition} \label{def:linkcongruence}
Let $\algA_1,\ldots,\algA_n$ be algebras and let $R \leq \prod_{i=1}^n \alg{A}_i$ be a relation with the parallelogram property. Then we define $\theta_i$ to be the relation:
\[\theta_i\ := \{(x,y) \in A_i^2 \mid  \exists z_1,\ldots, z_n R(z_1,\ldots,z_{i-1}, x, z_{i+1}, \ldots,z_n) \land R(z_1,\ldots,z_{i-1}, y, z_{i+1}, \ldots,z_n)\}.\]
\end{definition}

We call $\theta_i$ the \emph{$i$-th coordinate kernel} of $R$. We remark that this is sometimes also referred to as $i$-th \emph{link/linkedness congruence of $R$}~\cite{Zeb}.

\begin{lemma} \label{lemma:reduced}
Let $\algA_1,\ldots,\algA_n$ be algebras and let $R \leq \prod_{i=1}^n \alg{A}_i$ be a relation with the parallelogram property. Then the following holds: 

\begin{enumerate}
    \item $\theta_i$ is a congruence of $\pr_i(R) \leq \alg{A}_i$ for every $i=1,\ldots,n$.
    \item Let $\algA_i' = \pr_i(R)/\theta_i$, and $\nu \colon \prod_{i = 1}^n \pr_i(R) \to \prod_{i = 1}^n \alg{A}_i'$ be defined by $\nu(x_1,\ldots,x_n) = (x_1/\theta_1,\ldots,x_n/\theta_n)$. Then $R = \nu^{-1}(\nu (R))$.
\end{enumerate}
\end{lemma}

\begin{proof}
    We first prove (1). Without loss of generality, let $i=1$. Reflexivity, symmetry, and the fact that $\theta_1 \leq \alg A_1^2$ follow directly from the definition of $\theta_1$. For transitivity assume that $\theta_1(a_1,b_1), \theta_1(b_1,c_1)$. Then, by the definition of $\theta_1$, there exist $\boldsymbol{u}, \boldsymbol{v} \in \prod_{i = 2}^nA_i$ such that $R(x, \boldsymbol{u}), R(y, \boldsymbol{u}), R(y, \boldsymbol{v}), R(z, \boldsymbol{v})$. By the parallelogram property applied to the first 3 tuples, $R(x, \boldsymbol{v})$ holds, which implies $ \theta_1(x,z)$. 
Next we prove (2). By definition, $\boldsymbol{b} \in \nu^{-1}(\nu (R))$ if and only if there is an $\boldsymbol{a} \in R$ with $a_i \theta_i b_i$ for all $i=1,\ldots,n$. We are going to prove by induction on $i=0,1,\ldots,n$, that $R(b_1,\ldots,b_i,a_{i+1},\ldots,a_n)$. For $i=0$ this trivially holds. For an induction step, we can assume that $R(b_1,\ldots,b_{i-1},a_i,a_{i+1},\ldots,a_n)$. Then $a_i \theta_i b_i$ implies that there are elements $w_j \in A_j$ for $j\neq i$ such that $R(w_1,\ldots,w_{i-1},a_i,w_{i+1},\ldots,w_n)$ and $R(w_1,\ldots,w_{i-1},b_i,w_{i+1},\ldots,w_n)$. By the parallelogram property of $R$ this implies $R(b_1,\ldots,b_{i-1},b_i,a_{i+1},\ldots,a_n)$, which is what we wanted to prove. For $i=n$ we obtain $R(b_1,\ldots,b_n)$, which concludes the proof.
\end{proof}

\begin{definition}
Let $\algA_1,\ldots,\algA_n$ be algebras and let $R \leq \prod_{i=1}^n \alg{A}_i$ a relation with the parallelogram property, and $\nu$ be defined as in Lemma \ref{lemma:reduced}. Then we define the \emph{reduced representation} of $R$ as the image $\overline{R} = \nu(R)$.
\end{definition}

We are further going to need the following definition:

\begin{definition}
Let $\alg{A}$ be an algebra. An $n$-ary relation $R \leq \alg{A}^n$ is called \emph{critical} if it is \emph{$\land$-irreducible} in the lattice of subuniverses of $\alg{A}^n$ (i.e., it cannot be written as the intersection of strictly bigger relations $Q \leq \alg{A}^n$) and it has \emph{no dummy variables}, i.e., it depends on all of its inputs.
\end{definition}

Clearly $\Gamma$ is a relational basis of $\alg{A}$ if and only if all critical relations can be pp-defined from $\Gamma$. Thus, it is enough to consider only critical relations in our analysis. Theorem 2.5 in~\cite{KS-parallelogramterms} discusses the properties of critical relations with respect to algebras from congruence modular varieties. In the special case of Mal'cev algebras, the result directly implies the following (compare also with Theorem 2.3.10 in~\cite{Zeb}):

\begin{theorem}[Consequence of Theorem 2.5. in~\cite{KS-parallelogramterms}] \label{teor:critical}
Let $\alg{A}$ be a Mal'cev algebra, $n\geq 3$ and $R\leq \alg{A}^n$ be a critical relation. Further, let $\overline R \leq_{sd} \prod_{i=1}^n \alg{A}_i$ be the reduced representation of $R$. Then $\overline R$ is also critical (as a multi-sorted relation) and
\begin{enumerate}
\item $\alg{A}_i$ is subdirectly irreducible for every $i$,
\item the monolith $\mu_i$ of $\alg{A}_i$ is Abelian for every $i$,
\item for every pair of indices $i\neq j$, the quotient of $\pr_{i,j} \overline R$ by the centralizers $(0_{\alg{A}_i} : \mu_i) \times (0_{\alg{A}_j} : \mu_j)$ is the graph of an isomorphism $$\alg{A}_i/(0_{\alg{A}_i} : \mu_i) \cong \alg{A}_j/(0_{\alg{A}_j} : \mu_j).$$
\end{enumerate}
\end{theorem}

Note that Theorem~\ref{teor:critical} only applies to critical relations of arity 3 or higher. However, classifying the invariant relations of arity 2 is also an easy task: Let $R \leq \alg{A}^2$ and $\overline R \leq \alg{A}_1 \times \alg{A}_2$ be its reduced representation. Then $\overline R$ is an isomorphism between $\alg{A}_1$ and $\alg{A}_2$. Thus, for Mal'cev algebras, there is a one-to-one correspondence of binary invariant relations and the isomorphism between quotients of subalgebras of $\algA$.

Mal'cev algebras that only have binary critical relations can be characterized as follows.

\begin{theorem} \label{thm:majority}
Let $\algA$ be a Mal'cev algebra. The following are equivalent:
\begin{enumerate}
\item All critical relations $R\leq \algA^n$ are at most binary.
\item $\algA$ has a majority term.
\item there is no subdirectly irreducible $\algB \in \HS(\algA)$ with an Abelian monolith.
\end{enumerate}
\end{theorem}
\begin{proof}
By the Baker-Pixley Theorem~\cite{BakerPixley}, (1) is equivalent to (2) (even for non-Mal'cev algebras). For (2)$\Rightarrow$(3), note that the relation $T^\algB(\mu)$ for the monolith $\mu\neq 0_B$ of $\algB$ cannot be invariant under any majority operation $f$ (otherwise, evaluating $f$ at the tuples $(x,y,y,x), (y,y,y,y), (y,y,x,x) \in T^\algA(\mu)$ would imply $(y,y,y,x)\in T^\algA(\mu)$, so $x=d(y,y,y)$ for a Mal'cev term $d$ and all $x \mu y$ - contradiction). The implication (3)$\Rightarrow$(1) clearly follows from Theorem~\ref{teor:critical}.
\end{proof}

Varieties that have both a Mal'cev and majority operation are also called \emph{arithmetic} and are well-studied in the literature. In addition to Theorem~\ref{thm:majority}, several equivalent characterizations are known (e.g., by Pixley terms \cite[Theorem 4.70]{Bergman}, or as those Mal'cev varieties that are congruence neutral \cite[Theorem 7.50]{Bergman}).

We are further going to use the well-known fact that Abelian Mal'cev algebras of prime size are related by at most $4$-ary relations~\cite{BagyinskiDemetrovics82}. For better reference, we include a short proof sketch:

\begin{theorem} \label{thm:Zp}
Let $\alg{A}$ be an Abelian algebra with Mal'cev term $d$, and assume that $|A| = p$ is prime. A relational basis of $\alg{A}$ is given by $T^{\algA} \leq\alg{A}^4$, together with at most binary invariant relations.
\end{theorem}

\begin{proof}
By Definition~\ref{def:centralize}, $T^{\algA}$, the graph of the Mal'cev operation is invariant under $\algA$. Up to isomorphism, there is only one Abelian group of order $p$. Thus, by Lemma~\ref{lemma:centerwelldef} (3), we can assume without loss of generality that $A = \{0,1,\ldots,p-1\}$ and $T^{\algA} = \{(x,y,z,x-y+z) \mid x,y,z\in A\}$, where $+$ and $-$ denote addition and subtraction modulo $p$.

It is easy to see that $T^{\algA}$ pp-defines every $(n+1)$-ary relation of the form $\{ (\boldsymbol{x},z) \mid \ \sum_{i=1}^n \beta_ix_i = z \}$, for $\bm{\beta} \in \mathbb Z_p^n$ satisfying $\sum_{i=1}^n \beta_i = 1$. Indeed, these relations are exactly the graphs of the term operation of $\mathbf Z_p = (\{0,1,\ldots,p-1\}, x-y+z \bmod p)$.

Now, let $R \leq \alg{A}^n$ be an arbitrary invariant relation. As $R$ is invariant under $x-y+z$, it must be an affine subspace of the vector space $\mathbb{Z}_p^n$ (of some codimension $k$). Therefore, we can describe it by a system of linear equations $ \bm{\alpha}_j^t \boldsymbol{x} = c_j$,
 for $j = 1,\ldots,k$ (for linearly independent $\bm{\alpha}_j$). By the observation above, $T^\algA$ pp-defines every $(n+2)$-ary relation $Q_j = \{(\boldsymbol{x},y,z) \mid \bm{\alpha}_j^t \boldsymbol{x}  + (1-\sum_{i=1}^n \alpha_{i,j}) y = z \}$, for every $j$.

It is now straightforward to see that the binary relation $S_j(y,z)$, defined by $\exists \boldsymbol{x}( R(\boldsymbol{x}) \land Q_j(\boldsymbol{x},y,z))$, is equivalent to $z-(1-\sum_{i=1}^n \alpha_{i,j}) y = c_j$. Thus, $R$ can be pp-defined by the formula
\[\exists y_1,z_1,\ldots,y_k,z_k \bigwedge_{j = 1}^k Q_j(\boldsymbol{x},y_j,z_j) \land S_j(y_j,z_j).\]
This finishes the proof.
\end{proof}

It is a well-known fact that on a set of prime cardinality $p$ there is only one idempotent Abelian Mal'cev algebra up to term equivalence and isomorphisms (see e.g. \cite{BagyinskiDemetrovics82}; this can also straightforwardly be derived from Lemma \ref{lemma:centerwelldef}(3) and the fact that there is only one Abelian group of size $p$). In the rest of the paper, we will use the following notation for it:

\begin{definition}
Let $p$ be a prime. Then we denote by $\mathbf Z_p$ the algebra $\mathbf Z_p = (\{0,1,\ldots,p-1\}, x-y+z \bmod p)$, and by $\clo{Z}_p = \Clo(\mathbf Z_p)$, the clone of idempotent affine operations modulo $p$. Following the standard convention for relational structures, we further define the relational structure $\mathbb Z_p = (\{0,1,\ldots,p-1\},T^{\mathbf Z_p}, \{0\}, \ldots \{p-1\})$ (note that $\clo{Z}_p = \Pol(\mathbb Z_p)$).
\end{definition}

\section{Bulatov's relational bases} \label{sec:reprovingBulatov} 
In this section, we provide a self-contained proof that the at most 4-ary invariant relations of a three-element Mal'cev algebra $\alg{A}$ form a relational basis of $\alg{A}$. 
Moreover, we are going to show that each three-element Mal'cev algebra has a relational basis that is a subset of the very set of at most 4-ary relations given in \cite{BulatovImportante}.
We remark that the proof given in \cite{BulatovImportante}
is written in the language of tame congruence theory (TCT)~\cite{HobbyMcKenzie}, while we will use the language of commutator theory introduced in Section~\ref{sec:malcev}. For the interested reader, we remark that, for a finite Mal'cev algebra $\algA$, a pair of congruences $\alpha \prec \beta$ (i.e., that $\beta$ is a cover of $\alpha$) is of type $\tp{2}$ if $\beta/\alpha$ is Abelian in $\algA/\alpha$, and of type $\tp{3}$ otherwise. Note that no other TCT-types of congruence pairs appear in finitely generated Mal'cev varieties (cf. \cite[Theorem 9.14]{HobbyMcKenzie}).

We start our analysis of Mal'cev algebras on $\{0,1,2\}$ with two well-known facts about their congruence lattices. 

\begin{lemma} \label{lemma:alwaysSI}
Let $\alg{A}$ be a Mal'cev algebra over $\{0,1,2\}$. If $\mu$ is a nontrivial congruence, then it is unique. 
\end{lemma}
\begin{proof}
Seeking a contradiction, let us assume 
that both the equivalence relation $\mu_i$ with blocks $\{i\}$ and $\{0,1,2\} \setminus \{i\}$, and $\mu_j$ with blocks $\{j\}$ and $\{0,1,2\} \setminus \{j\}$ are congruences of $\alg{A}$, for $i\neq j$. Then formula $\exists y \mu_i(x_1,y) \land \mu_j(y,x_2)$ pp-defines the relation $\{0,1,2\}^2\setminus \{(i,j)\}$. This relation does not have the parallelogram property; hence $\algA$ cannot be Mal'cev.
\end{proof}

\begin{lemma} \label{lemma:centralizer}
Let $\alg{A}$ be a three-element Mal'cev algebra with an Abelian monolith $\mu$. Then $(0_{\alg{A}} : \mu) = \mu$.
\end{lemma}

\begin{proof}
The statement is trivially true if $\mathbf{A}$ is simple. So let us assume that $\mathbf{A}$ has a non-trivial congruence $\mu$. Then $\mu$ has one block of size 2 and one of size 1. As $\mu$ is Abelian, $(0_{\alg{A}} : \mu) \geq \mu$. If $(0_{\alg{A}} : \mu) = 1_{\alg{A}}$, then by Lemma~\ref{lemma:centerwelldef} (3), all $\mu$-classes $[a]_\mu$ need to have the same size, which is a contradiction. Thus, $(0_{\alg{A}} : \mu) = \mu$.
\end{proof}

We first study the at most binary invariant relations of $\alg{A}$. We recall from the previous section that these correspond exactly to the isomorphisms between quotients of subalgebras of $\algA$. In particular, this includes all isomorphisms between subalgebras of $\alg{A}$, which we are going to call \emph{partial isomorphisms}.

\begin{definition} \label{definition:newnotation}
We are going to use the following notation for binary relations on $\{0,1,2\}$:
\begin{enumerate}
\item By $\mu_i$ we denote the equivalence relation on $\{0,1,2\}$ with classes $\{i\}$ and $\{0,1,2\} \setminus \{i\}$.
\item By $\rho_i = \{0,1,2\}^2\setminus \mu_i$ we denote the complement of $\mu_i$.
\item The graphs of the transposition fixing $i$ are denoted by $\psi_i$. Moreover, $\psi_i'$ denotes its restriction to $\{0,1,2\}\setminus \{i\}$. For example,
$$\psi_2 = \begin{pmatrix}
    0 & 1 & 2 \\
    1 & 0 & 2
  \end{pmatrix},\ \psi_2' = \begin{pmatrix}
    0 & 1  \\
    1 & 0 
  \end{pmatrix}.$$
\item The graph of the cyclic permutation $(012)$ of $\{0,1,2\}$ is denoted by $\varphi$, and its restriction to $\{0,1,2\}\setminus \{i\}$ by $\varphi_i'$, i.e.,
$$\varphi = \begin{pmatrix}
    0 & 1 & 2 \\
    1 & 2 & 0
  \end{pmatrix} \text{ and } \varphi_0' = \begin{pmatrix}
    1 & 2 \\
    2 & 0
  \end{pmatrix},\ \varphi_1' = \begin{pmatrix}
    0 & 2 \\
    1 & 0
  \end{pmatrix},\ \varphi_2' = \begin{pmatrix}
    0 & 1 \\
    1 & 2
  \end{pmatrix}.$$
\item We further define the partial bijections
$$\varphi_3' = \begin{pmatrix}
    1 & 0  \\
    2 & 0 
  \end{pmatrix},\ \varphi_4' = \begin{pmatrix}
    0 & 1  \\
    2 & 1 
  \end{pmatrix},\ \varphi_5' = \begin{pmatrix}
    0 & 2  \\
    1 & 2 
  \end{pmatrix}.$$
\item Lastly, let us define the following two relations that describe maps from $\{0,1,2\}$ to $\{0,1\}$:
$$\tau_0 = \begin{pmatrix}
    0 & 1 & 2 \\
    0 & 0 & 1
  \end{pmatrix},\ \tau_1 = \begin{pmatrix}
    0 & 1 & 2 \\
    1 & 1 & 0
  \end{pmatrix}.$$
\end{enumerate}
\end{definition}

We remark that Definition~\ref{definition:newnotation} follows roughly the same notational conventions as~\cite{BulatovImportante}. However, we changed the way the relations are indexed for better readability. Note that every partial bijection that is not the identity restricted to some subset is equal to $\psi_i'$, $\varphi_i'$, or $(\varphi_i')^{-1}$ for some $i$.

We can narrow down the number of binary relations we need to consider by only constructing relational bases up to isomorphisms of $\algA$, i.e., up to renaming the elements of $\algA$. By Lemma~\ref{lemma:alwaysSI} and some easy observations, we then obtain the following.

\begin{lemma} \label{lemma:binary}
Let $\alg{B}$ be a Mal'cev algebra on a three-element set. Then $\alg{B}$ is isomorphic to an algebra $\algA$ on $\{0,1,2\}$ with monolith $\mu$, such that one of the following is true:
\begin{enumerate}
\item $\mu = 1_{\alg{A}}$ and every binary $R\leq \algA^2$ is pp-definable from the unary invariant relations together with
$$\{\varphi, \psi_2, \psi_0', \psi_1', \psi_2', \varphi_0', \ldots, \varphi_5' \} \cap \Inv(\algA),$$
\item $\mu = \mu_2$ and every binary $R\leq \algA^2$ is pp-definable from the unary invariant relations together with
$$\{\mu_2, \rho_2, \tau_0,\tau_1, \psi_2, \psi_2' \}\cap \Inv(\algA).$$
\end{enumerate}
\end{lemma}

\begin{proof}
Let $B = \{a,b,c\}$ be the universe of $\algB$. Let us first assume that $\algB$ is simple; we show that we fall within the scope of Case~(1). Recall that we can regard binary invariant relations as isomorphism between quotients of subalgebras of $\algB$. Since $\algB$ is simple and has only three elements, all binary invariant relations that are not pp-definable from unary invariant relations must be partial isomorphisms.
If $\algB$ has exactly one automorphism $(ab)$ of order $2$, then we obtain $\algA$ as the isomorphic copy under the map $a\mapsto 0,\ b\mapsto 1,\ c\mapsto 2$. In this case, $\psi_2$ is clearly invariant under $\algA$. Otherwise, we map $B$ to $\{0,1,2\}$ in an arbitrary way. In this case, the automorphism group of $\algA$ is either trivial or full; every corresponding binary relation is pp-definable from the list. Furthermore, the graph of every partial isomorphism of $\algA$ is equal to one of the given list or to one of its inverses. Thus, the given list of relations pp-defines every invariant $R\leq \algA^2$.

Next, we show that if $\mathbf{A}$ is not simple, then we fall within the scope of Case~(2). 
Let $\mu$ be the non-trivial congruence of $\algB$, and let us assume $a \mu b$. 
Furthermore, let $\alg{A}$ be the isomorphic copy of $\mathbf{B}$ under $a\mapsto 0$, $b\mapsto 1$, $c\mapsto 2$. 
Clearly $\mu_2$ is a non-trivial congruence of $\algA$, and by Lemma~\ref{lemma:alwaysSI} it is unique. 

If $R\leq \algA^2$ describes a homomorphism of $\algA$ to a proper subalgebra, it must have kernel $\mu_2$. If the image of the homomorphism is equal to $\{0,1\}$, then $R=\tau_0$ or $R=\tau_1$; in all other cases $R$ can be obtained by the restriction of $\mu_2$ or $\rho_2$ to some 2-element invariant subset (in the latter case, $\rho_2(x,y) \leftrightarrow \exists z_1,z_2 \mu_2(x,z_1) \land R(z_1,z_2)\land \mu_2(z_2,y)$ is pp-definable from $R$ and $\mu_2$, and therefore also invariant under $\algA$).

Next, observe that $\algA$ cannot have any other automorphism than $\psi_2$ (as composing it with the quotient map $\algA \to \algA/\mu_2$, would give rise to a homomorphism with other kernel than $\mu_2$).

Last, assume that $R \leq \algA^2$ is the graph of a partial isomorphism that is not equal to $\psi_2'$. If its domain and codomain intersect with both $\mu_2$-classes, then $R$ is a restriction of $\mu_2$ or $\rho_2$ to two 2-element subuniverses, and therefore pp-definable from them (and in the latter case, $\rho_2$ is also invariant under $\algA$). 

If $R$ is a partial isomorphism between $\{0,1\}$ and some other 2-element subuniverse, we can similarly see that it is pp-definable as a restriction of the inverses of $\tau_0$ or $\tau_1$.
\end{proof}

In the case where $\alg{A}$ does not have a critical relation of arity $3$ or greater, the relations in Lemma~\ref{lemma:binary} already form a relational basis for $\alg{A}$. This is exactly the case if $\HS(\alg{A})$ contains no algebras with an Abelian monolith; see Theorem~\ref{thm:majority}.

Next, we consider critical relations $R$ of arity $n\geq 3$.

\begin{lemma}\label{lemma:critical1}
Let $\algA$ be a Mal'cev algebra on a 3-element set with monolith $\mu$. Let $R\leq \alg{A}^n$ be a critical relation with $n\geq 3$, and let $\overline{R} \leq_{sd} \prod_{i=1}^n \alg{A}_i$ be its reduced representation. Then either
\begin{enumerate}
\item all $\alg{A}_i$ are Abelian, or
\item $\mu\neq 1_{\alg{A}}$, $\alg{A} = \alg{A}_i$ for all $i$, and $R$ is pp-definable from $T^\alg{A}(\mu)$ and at most binary invariant relations.
\end{enumerate}
\end{lemma}

\begin{proof}
First, assume that there is an index $i$ such that $\alg{A}_i$ is not Abelian. Since the monolith of $\algA_i$ is Abelian (Theorem~\ref{teor:critical} (2)), this is only possible if $\mu_i \neq 1_{A_i}$. Since $\algA_i$ has 3 congruences, this is only possible if $\alg{A}_i = \alg{A} = \pr_i(R)$. Thus, $\mu = \mu_i$ is Abelian, and hence, by Lemma~\ref{lemma:centralizer}, $(0:\mu) = \mu$. By Theorem~\ref{teor:critical} (3), every algebra $\algA_j$ must have a non-trivial quotient algebra of size 2. As $|A_j| \leq 3$, this must be the quotient by the mononlith $\mu_j \neq 1_{A_i}$. So, again Theorem~\ref{teor:critical} (2) implies that $\alg{A}_j=\alg{A} = \pr_j(R)$ for all $j \in [n]$. 

In the following, let us assume without loss of generality that $A = \{0,1,2\}$ and $\mu = \mu_2$. By Theorem~\ref{teor:critical} (3), $\pr_{i,j}(R)/\mu_2 \times  \mu_2$ is an automorphism of $\alg{A}/\mu_2$ for all $i \neq j$; hence, $\pr_{i,j}(R)$ is either equal to $\mu_2$ (if the automorphism is the identity on $\alg{A}/\mu_2$), or equal to $\rho_2$ (if the automorphism is the only transposition on $\alg{A}/\mu_2$).

We next show that the case $\pr_{i,j}(R) = \rho_2$ actually never appears. Seeking a contradiction, assume otherwise; then (up to reordering of coordinates), $R$ is a subset of $\{2\}^l \times \{0,1\}^{n-l} \cup \{0,1\}^{l} \times \{2\}^{n-l}$, and $\rho_2 = \pr_{l,l+1}(R)$ is invariant under $\alg{A}$. Then $R$ can be defined as  $\{\boldsymbol{x} \in \prod_{i=1}^n A_i \mid \pr_{1,\ldots,l}(R)(x_1,\ldots,x_l) \land \pr_{l+1,\ldots,n}(R)(x_{l+1},\ldots,x_n) \land \rho_2(x_l,x_{l+1})\}$; hence, $R$ is not critical - a contradiction.

It follows that $R$ is the disjoint union of the constant tuple $(2,\ldots,2)$ and $R\cap \{0,1\}^n$. By item (3) in Lemma~\ref{lemma:centerwelldef}, $d(x,y,z) = x-y+z\mod 2$ on $\{0,1\}$. Thus, the restriction of $\{0,1\}^n$ to $R$ needs to be equal to an affine subspace of the vector space $\mathbb{Z}_2^n$. In the same fashion as in the proof of Theorem~\ref{thm:Zp}, we can then show that $T^\alg{A}(\mu)$, together with all binary relations, pp-defines $R$.
\end{proof}

By Lemma~\ref{lemma:critical1} either all algebras in the reduced representation of a $n$-ary critical relation $R$ (with $n\geq 3$) are Abelian - Case~(1) or there exits a pp-definition of $R$ via standard relations - Case~(2).
One might think that  a pp-definition of $R$ from the relations $T^{\algA_i}$ and the binary relations can also be found in Case~(1) by applying Theorem~\ref{thm:Zp}. 
However, the Abelian algebras $\algA_i$ do not need to be equal or even isomorphic. We deal with the resulting technicalities in the following lemma:

\begin{lemma} \label{lemma:critical2}
Let $\algA$ be a Mal'cev algebra on a 3-element set with a monolith $\mu$. Let $R\leq \alg{A}^n$ be a critical relation with $n\geq 3$, and let $\overline{R} \leq_{sd} \prod_{i=1}^n \alg{A}_i$ be its reduced representation. Assume that all $\algA_i$ are Abelian. Then, one of the following holds:
\begin{enumerate}
\item $\algA$ is Abelian and $R$ is pp-definable from $T^\algA$ and at most binary relations, \label{itm:case1}
\item there is a two-element Abelian $\alg{B} \leq \alg{A}$, and $R$ is pp-definable from $T^\algB$ and at most binary relations, \label{itm:case2}
\item $\alg{A}/\mu$ is Abelian and $R$ is pp-definable from $\{(x_1,x_2,x_3,x_4) \in A^4 \mid (x_1/\mu,x_2/\mu,x_3/\mu,x_4/\mu) \in T^{\alg{A}/\mu}\}$ and at most binary relations, \label{itm:case3}
\item $\alg{A}/\mu$ is Abelian, the subalgebra $\{a,b\} \leq \alg{A}$ is Abelian and a block of $\mu$, and $R$ is pp-definable from binary relations 
and 
\[S_{ab} = \begin{pmatrix}
a & a & b & b & a & b\\
a & a & b & b & b & a\\
a & b & a & b & c & c
\end{pmatrix}.\] \label{itm:case4}
\end{enumerate}
\end{lemma}

\begin{proof}
We split the proof into cases: \\
\textbf{Case 1}: \textit{There exists $i$ such that $\alg{A}_i = \alg{A}$}:
A three-element Abelian Mal'cev algebra cannot have any subalgebras of size $2$ (this is e.g. a consequence of Lemma \ref{lemma:centerwelldef} (3)). Therefore 
$\alg{A}_j = \alg{A}$ holds for all $j\leq n$. Theorem~\ref{thm:Zp} for $p=3$
implies that condition~\eqref{itm:case1} is satisfied. 

Note that in all the remaining cases we have $|A_i|=2$ for all $i\leq n$, and $(R,d)$ is isomorphic to an affine subspace of $\mathbf{Z}_2^n$.

\textbf{Case 2}: \textit{Each $\alg{A}_i$ is a two element subalgebra of $\alg{A}$}:  
Let $\alg{B}$ be a two element subalgebra of $\alg{A}$ with $\alg{A}_r = \alg{B}$ for some $r\leq n$.
If each $\alg{A}_i$ is equal to $\alg{B}$, condition~\eqref{itm:case2} is satisfied by Theorem~\ref{thm:Zp}.
If there is $l\leq n$ such that $A_l\neq A_r$,
then, each $\alg{A}_i$ contains a 1-element subalgebra $c_i$, obtained by intersecting $\alg{A}_i$ with $\alg{A}_r$ or $\alg{A}_l$.
Let us assume that the affine space given by $R$ has dimension $k<n$. 
By setting $k-1$ variables of $R$ to invariant constants $c_i$ and 
then projecting the relation to some remaining coordinates $\{j,k\}$, 
we obtain an isomorphism between $\alg{A}_j$ and $\alg{A}_k$. 
Hence, all domains $\alg{A}_i$ are isomorphic.  
Up to composition with the isomorphisms, we can reduce to the case in which all domains are equal to the same subalgebra $\alg{B}$; this implies that condition~\eqref{itm:case2} is satisfied.

In the remaining cases we assume without loss of generality that $\mu = \mu_2$ an that there is $k\geq 1$
with $\alg{A}_1 = \cdots = \alg{A}_k = \alg{A}/\mu$, and $\alg{A}_j \leq \alg{A}$ for $j = k+1,\ldots,n$. 

\textbf{Case 3}: \textit{There exists $\alg{B}\leq\alg{A}$ and $h\colon \alg{A}\to \alg{B}$ with kernel $\mu$}: Then $R$ pp-defines the relation $R' = \{\boldsymbol{x} \in A^n \mid \exists \boldsymbol{y} R(y_1,\ldots,y_l,x_{l+1},\ldots,x_n) \land \bigwedge_{i=1}^l (x_i = h(y_i))\}$ which can be pp-defined from $T^{\algB}$ following the argument given for the previous case.
Since $R$ can be defined by the formula $\exists z_1,\ldots,z_l R'(z_1,\ldots,z_l,x_{l+1},\ldots,x_n) \land \bigwedge_{i=1}^l (z_i = h(x_i))$, also $R$ can be pp-defined by $T^{\alg B}$ and binary relations. Thus condition~\eqref{itm:case2} is satisfied. 

\textbf{Case 4}: \textit{There is no $\alg{B}\leq\alg{A}$ isomorphic to $\algA/\mu$}:
In particular this implies that $A_j = \{0,1\}$, for every $j>k$,
as otherwise the the second isomorphism theorem
would imply that $\algA_j \cong \algA/\mu$. 
We further split this case into subcases:

\textbf{Case 4.1}: \textit{$k=n$}:
Then $\overline R \leq (\alg{A}/\mu)^n$ can be pp-defined by $T^{\alg{A}/\mu}$ and binary relations, by Theorem~\ref{thm:Zp}. Hence condition~\eqref{itm:case3} is satisfied. 

\textbf{Case 4.2}: \textit{$A_i = B = \{0,1\}$ for all $i>k$}:
The reducts $(B,d) = (\{0,1\},x-y+z)$ and $(A/\mu,d^{\alg{A}/\mu})$ are isomorphic by Lemma~\ref{lemma:centerwelldef}(3) (and the fact that there is only one Abelian group of order 2).
Let $h \colon \{0,1,2\} \to \{0,1\}$ be a homomorphism from $(A,d)$ to $(B,d)$; without loss of generality let us assume $h(2) =1$ and $h(0)=h(1)=0$, otherwise we rename the elements of $A$ accordingly. Then $R$ can be represented by a system of linear equations of the form $\sum_{i=1}^k \alpha_i h(x_i) + \sum_{i=k+1}^n \beta_i x_i = c$.

We know that all relations given by equations $\{(\boldsymbol{x},z) \mid \sum_{i=1}^n \beta_i x_i = z\}$ with $\sum_{i=1}^n \beta_i = 1$ are pp-definable from $T^\algB$. Moreover, all relations of the form $\{(h(x_1),\ldots,h(x_n),h(z)) \mid \sum_{i=1}^n \alpha_i h(x_i) = h(z) \}$ are pp-definable from $T^{\alg{A}/\mu}$ and $0$ within $\algA/\mu$. The pre-image of these relations under $h$ are $\{(x_1,x_2,x_3,x_4) \in A^4 \mid (x_1/\mu,x_2/\mu,x_3/\mu,x_4/\mu) \in T^{\alg{A}/\mu}\}$ and $B = \{0,1\} \leq \algA$.
 
By applying the same substitution-trick as in the proof of Theorem~\ref{thm:Zp} (however separately for the variables $x_1,\ldots,x_k$ and the variables $x_{k+1},\ldots,x_n$), we get that $R$ is pp-definable from $T^\algB$, $B$, $\{(x_1,x_2,x_3,x_4) \in A^4 \mid (x_1/\mu,x_2/\mu,x_3/\mu,x_4/\mu) \in T^{\alg{A}/\mu}\}$, and ternary relations given by equations $\alpha h(x) + \beta_1y +\beta_2z = c$.

If one of the coefficients $\alpha$, $\beta_1$ or $\beta_2$ is $0$, such a ternary relation can clearly be pp-defined from an invariant binary relation. It is also easy to see that $h(x) +y+z = 0 \Leftrightarrow S_{01}(y,z,x)$. The only case left, i.e., the ternary relation given by $h(x) + y+z= 1$ is pp-definable from $S_{01}$ and $\psi_2'(y,z) \Leftrightarrow \exists x \in \{0,1\} \colon h(x)+ y + z = 1$. Thus $\alg{A}$ has a relational basis consisting of $T^\algB$, $\{(x_1,x_2,x_3,x_4) \in A^4 \mid (x_1/\mu,x_2/\mu,x_3/\mu,x_4/\mu) \in T^{\alg{A}/\mu}\}$, $S_{01}$ and all binary relations. It is an easy exercise to show that $S_{01}$ pp-defines $T^\algB$ and $\{(x_1,x_2,x_3,x_4) \in A^4 \mid (x_1/\mu,x_2/\mu,x_3/\mu,x_4/\mu) \in T^{\alg{A}/\mu}\}$, hence condition~\eqref{itm:case4} is satisfied.
\end{proof}

As a direct consequence of Theorem~\ref{teor:critical}, Lemma~\ref{lemma:critical1}, 
and Lemma~\ref{lemma:critical2}, we obtain that any Mal'cev algebra on a 3-element set 
has an at most 4-ary relational basis, given by the binary relations in Lemma~\ref{lemma:binary}, 
the witnesses of the Abelianness of $\algA$, its subquotients, and possibly a ternary relation $S_{ab}$ 
in the case that there is no $\alg{B}\leq\alg{A}$ isomorphic to $\algA/\mu$. 
In order to better restrict the set of at most 4-ary relations among which a relational basis can be chosen, we introduce the following notation, 
which again slightly differs from the one from \cite{BulatovImportante}:

\begin{definition}
By $\minority(x,y,z)$ we denote the \emph{minority} operation on a 2-element set, i.e., the operation defined by $\minority(x,y,y) = \minority(y,x,y) = \minority(y,y,x) = x$. On the set $\{0,1,2\}$, we define the relations:
\begin{itemize}
\item $T = \{(x,y,z, x-y+z \bmod 3) \mid x,y,z \in \{0,1,2\} \}$,
\item $T_{i}' = \{(x,y,z,\minority(x,y,z)) \mid x,y,z \neq i \}$, for every $i\in \{0,1,2\}$,
\item $T_{i} = T_i' \cup \{(i,i,i,i)\}$, for every $i\in \{0,1,2\}$,
\item $T^{\mu} = \{(x,y,z,u) \mid u/{\mu} = \minority(x/{\mu},y/{\mu},z/{\mu}) \}$, for every non-trivial equivalence relation $\mu$.
\end{itemize}
\end{definition}

Bulatov's description of all Mal'cev clones on a three-element set can be reformulated as follows:

\begin{theorem} \label{thm:bulatovreproved}
Let $\algB$ be a Mal'cev algebra on a three-element set. Then $\algB$ is isomorphic to an algebra $\algA$ with universe $\{0,1,2\}$ and monolith $\mu$, such that either
\begin{itemize}
\item $\mu = 1_A$, and $\algA$ is finitely related by its subuniverses together with
\begin{equation}\label{eq:standardrel}\{\varphi, \psi_2, \psi_0', \psi_1', \psi_2', \varphi_0', \ldots, \varphi_5', T,T_0',T_1',T_2'\}\cap \Inv(\algA), \end{equation}
\item or $\mu = \mu_2$, and $\algA$ is finitely related by its subuniverses together with
\begin{equation}\label{eq:standardrel2}\{\mu_2, \rho_2, \tau_0,\tau_1, \psi_2, \psi_2' ,T_2,T_2',T^{\mu_2}, S_{01} \} \cap \Inv(\algA). \end{equation}
\end{itemize}
\end{theorem}

\begin{proof}
We define the isomorphic copy $\algA$ of $\algB$ as in the proof of Lemma~\ref{lemma:binary}. Then, the theorem follows straightforwardly from Lemma~\ref{lemma:critical1}, Lemma~\ref{lemma:critical2}, and Lemma~\ref{lemma:binary}. If $\mu = \mu_2$, note that we do not need to include $T_0'$ or $T_1'$ in the relational basis, since if there is an Abelian subalgebra $\algB \leq \algA$ other than $\{0,1\}$, then $T^\algB$ is equal to the restriction of $T^{\mu_2}$ to $B$, and therefore pp-definable from the given basis.
\end{proof}

In the following, we will also refer to the relations in \eqref{eq:standardrel} or \eqref{eq:standardrel2} as the \emph{standard relations} or \emph{standard relational basis} of $\algA$.

In~\cite{BulatovImportante}, these relations were used (with the help of computer assistance) to obtain the full description of all 1129 Mal'cev clones on $\{0,1,2\}$ up to inclusion. We do not provide Bulatov's full classification here, but only discuss its minimal elements.

For each $i$ in $\{0,1,2\}$, let us define the Mal'cev operations 
 \begin{equation}\label{equation:thed_ifunctions}
 d_i(x,y,z)=\begin{cases}
   \minority (x,y,z) & \text{if } |\{x, y, z\}|\leq 2,\\
   i & \text{if } |\{x, y, z\}|=3.
   \end{cases}
\end{equation}
Moreover, we let 
\begin{equation}\label{equation:minority_plu_porj}
 g(x,y,z)=\begin{cases}
   \minority (x,y,z) & \text{if } |\{x, y, z\}|\leq 2,\\
   x & \text{if } |\{x, y, z\}|=3
   \end{cases}
\end{equation}

Then, the following holds:

\begin{theorem}[{\cite[Corollary~3.5]{BulatovImportante}}] \label{thm:minimalclones}
There are exactly five Mal'cev clones on $\{0,1,2\}$ that are minimal with respect to inclusion, namely:
\begin{itemize}
\item $\mathcal Z_3 = \Pol(\mathbb Z_3) = \Pol(T,\{0\},\{1\},\{2\})$,
\item $\Clo(g)$ -- all standard relations in \eqref{eq:standardrel} except for $T$ form an relational basis,
\item $\Clo(d_i)$, for $i=0,1,2$. For $i=2$, all standard relations in \eqref{eq:standardrel2} form an relational basis.
\end{itemize}
\end{theorem}

Interestingly, the aforementioned clones also turn out to be minimal Taylor clones over $\{0,1,2\}$. The latter claim follows from a recent classification of all minimal Taylor clones over $\{0,1,2\}$, which can be found in~\cite{Zeb,MinTaylorOver3}. The question of whether it is true that every minimal Mal'cev clone is a minimal Taylor clone also for universes with more than three elements still remains open.

\section{Boolean and self-dual Mal'cev clones up to minion homomorphism} \label{sec:selfdual}

In this section, we present, for the reader's and our convenience, classifications of all Mal'cev clones on $\{0,1\}$ and all self-dual Mal'cev clones on $\{0,1,2\}$ up to minion homomorphism. These results follow from only considering the Mal'cev clones in the classification in \cite{albert} and \cite{BodirskyVucajZhuk} respectively. Let us recall from the introduction that, in our paper, $\struct{C}_2$ denotes the relational structure $\struct{C}_2 = (\{0,1\},\neq,\{0\}, \{1\})$. We remark that by Proposition \ref{prop:idempotent} $\struct{C}_2$ is pp-interconstructible with $(\{0,1\},\neq)$, which in the literature is also occasionally called $\struct{C}_2$.

\begin{theorem}[Consequence of Theorem 3.22. in~\cite{albert}] \label{thm:2element}
Let $\mathcal C$ be a Mal'cev clone on a two-element set. Then $\mathcal C$ is minor-equivalent to exactly one of the following clones:
\begin{itemize}
\item $\mathcal T$, the clone on a 1-element set  \hfill (denoted by $[0]$ in~\cite{albert}),
\item $\mathcal I_2 = \Pol(\{0\},\{1\})$  \hfill (denoted by $[m,q]$ in~\cite{albert}),
\item $\mathcal C_2 = \Pol(\mathbb C_2)$  \hfill (denoted by $[d_3,m]$ in~\cite{albert}),
\item $\mathcal Z_2 = \Pol(\mathbb Z_2)$  \hfill (denoted by $[m]$ in~\cite{albert}).
\end{itemize}
Moreover, $\mathcal Z_2 \isaminorof \mathcal C_2 \isaminorof \mathcal I_2 \isaminorof \mathcal T$.
\end{theorem}

\begin{theorem}[Consequence of Theorem 5.8 in~\cite{BodirskyVucajZhuk}] \label{thm:selfdual}
Let $\mathcal C$ be an idempotent Mal'cev clone on the three-element set $\{0,1,2\}$. If $\mathcal C$ preserves $\varphi$, then it is minor-equivalent to exactly one of the following clones:
\begin{itemize}
\item $\mathcal Z_3 = \Pol(\{0,1,2\};T, \{0\},\{1\},\{2\} )$  \hfill (denoted by $\mathbf {L}_3$ in~\cite{BodirskyVucajZhuk}),
\item $\mathcal L_2 = \Pol(\{0,1,2\};\varphi,\psi_2,T_2)$  \hfill (denoted by $\mathbf {TL}_2$ in~\cite{BodirskyVucajZhuk}),
\item $\mathcal D = \Pol(\{0,1,2\};\varphi, \psi_2', \{0\},\{1\},\{2\})$  \hfill (denoted by $\mathbf D$ in~\cite{BodirskyVucajZhuk}),
\item $\mathcal C_3 = \Pol(\{0,1,2\};\varphi, \{0\},\{1\},\{2\})$ \hfill (denoted by $\mathbf C_3$ in~\cite{BodirskyVucajZhuk}).
\end{itemize}
By inclusion, $\mathcal Z_3 \isaminorof \mathcal C_3$, and $\mathcal L_2 \isaminorof \mathcal D \isaminorof \mathcal C_3$ hold. By the minor conditions in Figure~\ref{fig:selfdual}, no other comparison holds.
\end{theorem}

We also observe that no clone that preserves $\varphi$ has a 3-cyclic term. Since all Taylor clones on a 2-element set have a 3-cyclic term, none of the clones in Theorem~\ref{thm:selfdual} is minor-equivalent to a clone from Theorem~\ref{thm:2element}. This already gives rise to 8 distinct elements in $\Pthree$ with a quasi Mal'cev operation. 

\begin{figure}[h]
\centering
\begin{tabular}{r|llll}
& $\clo{Z}_3 \not \models$ & $\clo{L}_2 \not \models$ & $\clo{D} \not \models$ \\
\hline
$\clo{Z}_3 \models$ & & 2-cyc & 2-cyc \\
$\clo{L}_2 \models$ & min & & \\
$\clo{D} \models$ & min & maj & \\
$\clo{C}_3 \models$ & min & maj & 2-cyc
\end{tabular}
\caption{Minor conditions separating clones of self-dual operations. Here, \emph{min} and \emph{maj} denote the minor conditions \emph{quasi minority} and \emph{quasi majority}, respectively.}
\label{fig:selfdual}
\end{figure}

\section{Algebras with a majority term} \label{sec:SDmeet}

We next focus on Mal'cev algebras that have a majority term. 
Recall that, by Theorem~\ref{thm:majority}, 
this implies that their clone of term operations has a binary relational basis. 
If the clone of term operations of a Mal'cev algebra preserves $\varphi$, 
then it belongs to one of two $\minorequivalentto$-classes by Theorem~\ref{thm:selfdual}, 
namely $\clo{D}/\minorequivalentto$ and $\clo{C}_3/\minorequivalentto$. 
Thus, in the rest of this section, we only consider  those clones that do not preserve $\varphi$. 
We distinguish two cases: simple algebras and those with a non-trivial monolith $\mu_2$.

\subsection{Simple algebras}

By Theorem~\ref{thm:bulatovreproved} (1), every simple Mal'cev clone $\clo{C}$ with majority that does not preserve $\varphi$ must be the polymorphism clone of its unary invariant relations, together with a subset of $\{\psi_2,\psi'_0,\psi'_1,\psi'_2,\{\varphi'_j\}_{j\in[6]}\}$. In other words
$\clo{C} \supseteq \Pol(\struct{C}^{*}_2)$, for 
\[\struct{C}^{*}_2=\Big(\{0,1,2\};\psi_2,\psi'_0,\psi'_1,\psi'_2,\{\varphi'_j\}_{j\in[6]},\{0,1\},\{0,2\},\{1,2\},\{i\}_{i\in\{0,1,2\}}\Big).\]

We show that $\struct{C}_2 = (\{0,1\};\neq,\{0\},\{1\})$ pp-constructs $\struct{C}_2^{*}$, which implies that $\mathcal C_2 = \Pol(\struct{C}_2) \isaminorof \mathcal C$.

\begin{lemma}\label{lem:C2Collapse}
The relational structure $\struct{C}_2$ pp-constructs $\struct{C}^{*}_2$.
\end{lemma}
\begin{proof}

We first observe that the structure
$$\struct{C}^{**}_2 = \Big(\{0,1,2\};\psi_2,\varphi_2'\Big)$$
pp-defines all relations of $\struct{C}^{*}_2$. For this, first note that we can pp-define $\{0,1\}$ and $\{1,2\}$ by projecting $\varphi_2'$ to its coordinates. By restricting $\psi_2$ to $\{0,1\}$, respectively $\{1,2\}$, we get $\psi_2'$, respectively $\varphi_5'$. Thus, we pp-defined graphs of bijections between all two element subsets. We can obtain any other partial bijection as a compositions of these partial bijections and the transposition $\psi_2'$. Moreover, all 2-element subsets are clearly pp-definable (by projecting these maps to their coordinates), as well as all singleton relations $\{i\}$ as intersections of 2-element subsets.

Since $\struct{C}^{*}_2$ is pp-definable from $\struct{C}^{**}_2$, it is enough to prove that $\struct{C}_2$ pp-constructs $\struct{C}^{**}_2$. For this, we need to construct a pp-power $\struct{S}$ of $\struct{C}_2$ and show that there exists a homomorphism from $\struct{S}$ to $\struct{C}^{**}_2$.  We define $\struct{S} = (\{0,1\}^2;\Psi_2,\Phi'_2)$ by
\begin{align*}
\Psi_2(\boldsymbol{a},\boldsymbol{b}) &\Leftrightarrow a_1 \neq b_2 \wedge a_2 \neq b_1,\\
\Phi_2'(\boldsymbol{a},\boldsymbol{b}) &\Leftrightarrow a_1=a_2\land a_2\neq b_2\land b_1=1.
\end{align*}

Further, we define the maps
\begin{align*}\label{equation:mappahom1}
 \iota(x)=\begin{cases}
   (0,0) & \text{if } x = 0,\\
   (1,1) & \text{if } x = 1,\\
   (1,0) & \text{if } x = 2.
   \end{cases}
 &&\nu(\boldsymbol{x})=\begin{cases}
   0 & \text{if } \boldsymbol{x} = (0,0),\\
   1 & \text{if } \boldsymbol{x} = (1,1),\\
   2 & \text{if } x_1 \neq x_2.
   \end{cases}
\end{align*}
It is straightforward to see that $\iota$ is a homomorphism from $\struct{C}^{**}_2$ to $\struct{S}$, as 
\begin{align*}
\iota(\psi_2) &= \{ ((0,0),(1,1)), ((1,1),(0,0)), ((1,0),(1,0))  \} \subseteq \Psi_2,\\
\iota(\varphi_2') &= \{ ((0,0),(1,1)), ((1,1),(1,0))  \} \subseteq \Phi_2',
\end{align*}

so, it is only left to show that $\nu$ is a homomorphism from $\struct{S}$ to $\struct{C}^{**}_2$.

Suppose $(\boldsymbol{a},\boldsymbol{b})\in\Psi_2$. By the definition of $\Psi_2$, one of two cases holds: either $a_1=a_2 \neq b_1 = b_2$ or $a_1=b_1 \neq a_2 = b_2$. In the first case, clearly $\nu(\boldsymbol{a}),\ \nu(\boldsymbol{b})\in\{0,1\}$ and $\nu(\boldsymbol{a})\neq\nu(\boldsymbol{b})$ hold; thus $(\nu(\boldsymbol{a}),\nu(\boldsymbol{b}))\in\psi_2$. In the second case $(\nu(\boldsymbol{a}),\nu(\boldsymbol{b}))=(2,2)\in\psi_2$.

Next, suppose $(\boldsymbol{a},\boldsymbol{b})\in\Phi_2'$. By definition, this holds if and only if $a_1= a_2 \neq b_2$. If $a_1=1$, then $\boldsymbol{a} = (1,1)$, $\boldsymbol{b} = (1,0)$, and thus $(\nu(\boldsymbol{a}),\ \nu(\boldsymbol{b}))= (1,2) \in \varphi_2'$. If $a_1=0$, then $\boldsymbol{a} = (0,0)$, $\boldsymbol{b} = (1,1)$, so $\nu(\boldsymbol{a}),\ \nu(\boldsymbol{b})= (0,1) \in \varphi_2'$. 

It follows that $\nu$ is a homomorphism, which finishes the proof.
\end{proof}

As a direct consequence, we obtain the following lemma:

\begin{lemma} \label{lemma:simple}
Let $\mathcal C$ be an idempotent Mal'cev clone on $\{0,1,2\}$ that does not preserve  $\mu_1,\mu_2,\mu_3$ and $\varphi$. Then $\mathcal C \minorequivalentto \mathcal C_2$ or $\mathcal C \minorequivalentto \mathcal I_2$.
\end{lemma}

\begin{proof}
This follows directly from Lemma~\ref{lem:C2Collapse} and the fact that in $\mathfrak{P}_3$ (and even in $\Pfin$) there are no elements between $\overline{\clo{C}_2}$ and $\overline{\clo{I}_2}$~\cite{VucajZhuk, meyerSubmaximal}.
\end{proof}

\subsection{Algebras with a non-trivial monolith}
We proceed by examining the case of Mal’cev algebras with a majority term whose associated algebra is not
simple. 
Without loss of generality, we will always assume that the unique non-trivial congruence is $\mu_2$. 
By Theorem~\ref{thm:bulatovreproved} (2), the only other possible invariant standard relations 
are the unary ones and $\psi_2$, $\psi_2'$, $\rho_2$, $\varphi$, $\tau_1$, and $\tau_2$. 
We are going to prove that, up to minor-equivalence, there are exactly four clones with this property (Theorem~\ref{theorem:classesOfSDmeetMalcev}). 
The Lemmata~\ref{lemma:Hppconstruction},~\ref{lemma:extensionm1},~\ref{lem:forTheLastCollapse}.
provide the necessary pp-constructions.

First, we define the structure 
\[\mathbb H = (\{0,1,2\}; \mu_2, \tau_0,\tau_1, \psi_2', \rho_2, \{0,2\}, \{1,2\}, \{0\},\{1\},\{2\}).
\] 
Note that $\mathbb H$ does not pp-define $\psi_2$, and, in fact, is the maximal subset of binary standard relations that pp-defines $\mu_2$, but does not pp-define $\psi_2$.

\begin{lemma} \label{lemma:Hppconstruction}
$\mathbb C_2$ pp-constructs $\mathbb H$ and vice-versa. 
\end{lemma}

\begin{proof}
It is straightforward to see that $\mathbb H$ pp-constructs $\mathbb C_2$, since $\{0\}$, $\{1\}$ and $\psi_2'$ (the inequality restricted to $\{0,1\}$) are relations of $\mathbb H$, and $\{0,1\}$ is pp-definable by the formula $\mu_2(x,0)$.

Next, we prove that $\mathbb C_2$ pp-constructs $\mathbb H$. We first observe that all relations in $\struct{H}$ are already pp-definable from $\mathbb H'= (\{0,1,2\}; \tau_0,\psi_2', \{0,2\}, \{1,2\})$: We obtain $\{0,1\}$ as projection of $\psi_2'$ to one of its coordinates. All singletons $\{i\}$ can be obtained as intersections of 2-element sets. We can further pp-define $\tau_1(x,y)$ by $\exists z (\tau_0(x,z) \land \psi_2'(z,x))$, $\rho_2(x,y)$ by $\exists z (\tau_0(x,z) \land \tau_1(y,z))$, and $\mu_2(x,y)$ by $\exists z (\rho_2(x,z) \land \rho_2(z,y))$.

In order to show that $\mathbb C_2$ pp-constructs $\mathbb H'$, let us define the following pp-power of $\mathbb C_2$:
$$\mathbb S = (\{0,1\}^2; \mathcal{T}_0,\Psi_2', U_{02}, U_{12})$$ where the relations are defined by
\begin{align*}
\boldsymbol{a} \in U_{02}& \Leftrightarrow a_1=0,\\
\boldsymbol{a} \in U_{12}& \Leftrightarrow a_1 = a_2,\\
(\boldsymbol{a},\boldsymbol{b}) \in \mathcal{T}_0& \Leftrightarrow
b_2 = 1 \wedge b_1 \neq a_2,\\
(\boldsymbol{a},\boldsymbol{b}) \in \Psi_2'& \Leftrightarrow a_2 = b_2 =1 \wedge  a_1 \neq b_1.
\end{align*}
By its definition, $\mathbb S$ is a pp-power of $\mathbb C_2$. It is straightforward to see that the map $h \colon \mathbb H \to \mathbb S$ defined by $h(0) = (0,1)$,\ $h(1) = (1,1)$,\ $h(2) = (0,0)$ is a homomorphism.

Next, let us define the map $g \colon \mathbb S \to \mathbb H$ defined by $$g(\boldsymbol{x}) = \begin{cases} 2 &\text{if } x_2 = 0,\\
1 &\text{if } \boldsymbol{x} = (1,1),\\
0 &\text{if } \boldsymbol{x} = (0,1).
\end{cases}$$
We are going to show that $g$ is a homomorphism. By definition, $g(U_{02}) = \{0,2\}$ and $g(U_{12}) = \{1,2\}$. For a pair $(\boldsymbol{a},\boldsymbol{b}) \in \mathcal T_0$, note that $g(\boldsymbol{b}) \in \{0,1\}$, since $b_2=1$. If $a_2 = 0$, then $b_1\neq a_2$ implies that  $\boldsymbol{b}=(1,1)$. Thus, $(g(\boldsymbol{a}),g(\boldsymbol{b})) = (2,1) \in \tau_0$. Similarly, if we assume $a_2 = 1$, then $\boldsymbol{b}=(0,1)$. In this case, $(g(\boldsymbol{a}),g(\boldsymbol{b}))$ is equal to $(0,1)$ or $(1,1)$, which are both elements of $\tau_0$. Thus, $g(\mathcal T_0) = \tau_0$. It is furthermore straightforward to show that 
$g(\Psi_2') = \{(0,1),(1,0)\} = \psi_2'$. 

In conclusion, $g$ is a homomorphism. Hence, $\mathbb H'$ is homomorphically equivalent to a pp-power of $\mathbb C_2$, which finishes the proof.
\end{proof}

Next, let us define the structures: 
\begin{align*}
\struct{M}_1 &= (\{0,1,2\};\psi_2,\mu_2,\{0\},\{1\},\{2\}),\\
\struct{M}''_1 &= (\{0,1,2\};\psi_2,\psi_2',\mu_2,\{0,1\},\{1,2\},\{0,2\},\{0\},\{1\},\{2\}).
\end{align*}

Note that by Theorem~\ref{thm:bulatovreproved},
each idempotent Mal'cev clone $\mathcal C$ that has a majority operation
and preserves $\mu_2$ and $\psi_2$, 
but not $\rho_2$, lies in the interval 
$\Pol(\struct{M}_1'') \subseteq \mathcal C \subseteq \Pol(\struct{M}_1)$. 
The following lemma shows that all the clones in this interval belong to the same equivalence class 
modulo $\minorequivalentto$. 

\begin{lemma} \label{lemma:extensionm1}
$\struct{M}_1$ pp-constructs $\struct{M}''_1$, and vice-versa.
\end{lemma}

\begin{proof}
By inclusion, $\struct{M}''_1$ pp-constructs $\struct{M}_1$, thus, we only need do show that $\struct{M}_1$ pp-constructs $\struct{M}''_1$.

Note that $\{0,1\}$ has a pp-definition in $\mathbb M_1$ by $x \in \{0,1\} \Leftrightarrow  \, (x,y) \in  \mu_2\wedge y \in \{0\}$; so also $\psi_2'$ also a pp-definition in $\mathbb M_1$ by $(x,y) \in \psi_2' \Leftrightarrow (x,y) \in\mu_2 \land x \in \{0,1\}$. Therefore is enough to show that $\struct{M}'_1$ pp-constructs $\struct{M}''_1$, where 
\[\struct{M}'_1 = (\{0,1,2\};\psi'_2, \psi_2, \mu_2, \{0,1\}, \{0\},\{1\},\{2\}).\]

We define the pp-power $\struct{S}$ of $\struct{M}'_1$ by
\[\struct{S}=(\{0,1,2\}^3;\Psi'_2,\Psi_2,\mathcal M_2,U_{01},U_{02},U_{12},U_0,U_1,U_2)\]
where $U_0 = \{(0,0,1)\}$, $U_1 = \{(1,0,1)\}$, $U_2 = \{(2,2,2)\}$, and
\begin{align*}
\boldsymbol{a}\in U_{01} &\Leftrightarrow\psi'_2(a_2,a_3)\wedge a_1\in\{0,1\},
    \\\boldsymbol{a} \in U_{02} &\Leftrightarrow \psi_2(a_2,a_3) \wedge a_1=a_2,
    \\ \boldsymbol{a} \in U_{12} &\Leftrightarrow \psi_2(a_2,a_3)\wedge a_1=a_3,\\
    \Psi'_2(\boldsymbol{a},\boldsymbol{b}) &\Leftrightarrow \psi'_2(a_1,b_1) \wedge \psi'_2(a_2,a_3) \wedge \psi'_2(a_2,b_3)\ \wedge \psi'_2(a_3,b_2)\wedge \psi'_2(b_2,b_3),
    \\\Psi_2(\boldsymbol{a},\boldsymbol{b}) &\Leftrightarrow \psi_2(a_1,b_1) \wedge \psi_2(a_2,a_3) \wedge \psi_2(a_2,b_3)\ \wedge \psi_2(a_3,b_2)\wedge \psi_2(b_2,b_3),\\
    \mathcal M_2(\boldsymbol{a},\boldsymbol{b}) &\Leftrightarrow \mu_2(a_1,b_1)\wedge \mu_2(a_1,a_2)\wedge \mu_2(b_1,b_2)\ \wedge \psi_2(a_2,a_3)\wedge \psi_2(b_2,b_3).
\end{align*}
We further define the maps:
\begin{align*}\label{eq:mapHom3}
 \delta(x)=\begin{cases}
   (0,0,1) & \text{if } x = 0,\\
   (1,0,1) & \text{if } x = 1,\\
   (2,2,2) & \text{if } x = 2.
   \end{cases}
 &&\nu(\boldsymbol{x})=\begin{cases}
   0 & \text{if } \boldsymbol{x} \in \{(0,0,1),(1,1,0)\},\\
   1 & \text{if } \boldsymbol{x} \in \{(0,1,0),(1,0,1)\},\\
   2 & \text{otherwise}.
   \end{cases}
\end{align*}
We first prove that $\delta$ is a homomorphism from $\struct{M}''_1$ to $\struct{S}$.

Suppose $(a,b)\in\psi'_2$, i.e., $a,b\in\{0,1\}$ and $a\neq b$. It follows that $(\delta(a),\delta(b)) = (\boldsymbol{a}, \boldsymbol{b})\in\{((0,0,1),(1,0,1)),((1,0,1),(0,0,1))\}$ and all the conjuncts of $\Psi'_2$ are satisfied, yielding $(\delta(a),\delta(b))\in\Psi'_2$.

Suppose $(a,b)\in\psi_2$, that is: either $a,b\in\{0,1\}$ and $a\neq b$ or $a=b=2$. In the first case $(\delta(a),\delta(b)) = (\boldsymbol{a}, \boldsymbol{b}) \in \psi_2$ by the above paragraph. If $a=b=2$, then $(\delta(a),\delta(b))=((2,2,2),(2,2,2)) \in \Psi_2$.

Suppose $(a,b)\in\mu_2$, that is: either $(a,b)\in\{0,1\}^2$ or $a=b=2$, and let $(\delta(a),\delta(b)) = (\boldsymbol{a}, \boldsymbol{b})$. Note that, by the definition of $\delta$, $\psi_2(a_2,a_3)\wedge\psi_2(b_2,b_3)$ always holds. It is straightforward to check that the first three conjuncts hold. Thus $(\delta(a),\delta(b))\in \mathcal M_2$. It is also straightforward to check that $\delta$ also preserves all the unary relations. Thus $\delta$ is a homomorphism.

Next, we show that $\nu$ is a homomorphism from $\struct{S}$ to $\struct{M}''_1$. First suppose that $(\boldsymbol{a},\boldsymbol{b})\in\Psi'_2$. By the definition of $\Psi'_2$, it follows that $\boldsymbol{a},\boldsymbol{b}\in\{(0,0,1),$ $(1,1,0),(0,1,0), (1,0,1)\}$. If $\boldsymbol{a}=(0,0,1)$, then $\boldsymbol{b}$ must be equal to $(1,0,1)$ and we have $(\nu(\boldsymbol{a}),\nu(\boldsymbol{b}))=(0,1)\in\psi'_2$. Similarly, if $\boldsymbol{a}=(1,1,0)$, the $\boldsymbol{b}$ must be equal to $(0,1,0)$ and we have $(\nu(\boldsymbol{a}),\nu(\boldsymbol{b}))=(0,1)\in\psi'_2$. The two remaining cases are covered by the fact that $\Psi'_2$ is symmetric.

Suppose $(\boldsymbol{a},\boldsymbol{b})\in\Psi_2$. If $\nu(\boldsymbol{a})$ belongs to $\{0,1\}$, then so does $\nu(\boldsymbol{b})$; by the definition of $\Psi_2$ and the fact that $\psi_2$ coincides with $\psi'_2$ when restricted to $\{0,1\}$, it follows that $(\nu(\boldsymbol{a}),\nu(\boldsymbol{b}))\in\psi_2$. If $\nu(\boldsymbol{a})=2$, then we distinguish two cases: either $a_i=2$, for some $i\in[3]$ or $\boldsymbol{a}\in\{0,1\}^3$ with $a_2=a_3$. In the first case, there exist $j\in[3]$ such that $b_j=2$; thus $\nu(\boldsymbol{b})=2$ and $(\nu(\boldsymbol{a}),\nu(\boldsymbol{b}))=(2,2)\in\psi_2$. For the second case we can see that there is no $(\boldsymbol{a},\boldsymbol{b})\in\Psi_2$ satisfying $a_2 = a_3$ with $\boldsymbol{a}\in\{0,1\}^3$. 

Suppose $(\boldsymbol{a},\boldsymbol{b})\in \mathcal M_2$ and that $(\nu(\boldsymbol{a}),\nu(\boldsymbol{b}))\notin\mu_2$; say $\nu(\boldsymbol{a})\in\{0,1\},\ \nu(\boldsymbol{b})=2$. That is $a_1,a_2,a_3\in\{0,1\}$ and $a_2\neq a_3$. By the definition of $\mathcal M_2$, it follows that $b_1,b_2,b_3\in\{0,1\}$ and $b_2\neq b_3$. Hence, $\nu(\boldsymbol{b})\in\{0,1\}$, a contradiction. Note that, since $\mathcal M_2$ is symmetric by definition, this is the only case we need to check.

As for the unary relations: suppose $\boldsymbol{a}\in U_{01}$, i.e., $\boldsymbol{a}\in\{0,1\}^3$ and $a_2\neq a_3$. It immediately follows from the definition of $\nu$ that $\nu(\boldsymbol{a})\in\{0,1\}$. Suppose $\boldsymbol{a}\in U_{02}$, i.e., $a_1=a_2$ and $\psi_2(a_2,a_3)$. If $a_1\in\{0,1\}$, then $\boldsymbol{a}$ is of the form $(a_1,a_1,1-a_1)$ and $\nu(\boldsymbol{a})=0$; if $a_1=2$, then $\boldsymbol{a}=(2,2,2)$ and $\nu(\boldsymbol{a})=2$. Finally, suppose $\boldsymbol{a}\in U_{12}$, i.e., $a_1=a_3$ and $\psi_2(a_2,a_3)$. If $a_1\in\{0,1\}$, then $\boldsymbol{a}$ is of the form $(a_1,1-a_1,a_1)$ and $\nu(\boldsymbol{a})=1$. Similarly as before, if $a_1=2$, then $\boldsymbol{a}=(2,2,2)$ and $\nu(\boldsymbol{a})=2$. Clearly $\nu$ also preserves the singleton relations. Thus $\nu$ is a homomorphism.

In conclusion, we have shown that $\mathbb M_1''$ is homomorphically equivalent to the pp-power $\mathbb S$ of $\mathbb M_1'$, thus we are done.
\end{proof}

For the next lemma, we define the following relational structures:
\begin{align*}
    \struct{M}_0 &= (\{0,1,2\};\psi_2,\rho_2,\{0,1\},\{0\},\{1\},\{2\}),
    \\\struct{M}_0' &= (\{0,1,2\};\psi_2,\rho_2,\tau_0,\tau_1,\{0,1\},\{0\},\{1\},\{2\}).
\end{align*}
 Note that $\Pol(\struct{M}_0)$ is the biggest idempotent Mal'cev clone with a majority operation that preserves $\psi_2$, $\rho_2$ and hence also $\mu_2$ (as $\mu_2(x,y)$ is equivalent to $\exists z (\rho_2(x,z) \land \rho_2(z,y))$). Thus, the interval of clones between $\Pol(\struct{M}_0')$ and $\Pol(\struct{M}_0)$ contains all idempotent Mal'cev clones with majority that preserve $\mu_2$, $\psi_2$, and $\rho_2$. The following lemma shows that this interval collapses with respect to $\minorequivalentto$.

\begin{lemma}\label{lem:forTheLastCollapse}
The relational structure $\struct{M}_0$ pp-constructs $\struct{M}_0'$, and vice-versa.
\end{lemma}
\begin{proof}
By inclusion $\struct{M}_0'$ clearly pp-constructs $\struct{M}_0$. 

For the converse, we define a pp-power $\struct{S}$ of $\struct{M}_0$ by
\[\struct{S} = (\{0,1,2\}^3;\Psi_2,\mathcal R_2,\mathcal T_0,\mathcal T_1,U_{01},U_0,U_1,U_2),
\]
where $U_0 =  \{(0,0,2)\}, U_1 = \{(1,2,1)\}, U_2 = \{(2,0,1)\}$, and
\begin{align*}
    \Psi_2(\boldsymbol{a},\boldsymbol{b}) &\Leftrightarrow \psi_2(a_1,b_1)\wedge \psi_2(a_2,b_3)\wedge \psi_2(a_3,b_2),
    \\\mathcal R_2(\boldsymbol{a},\boldsymbol{b}) &\Leftrightarrow \rho_2(a_1,b_1),
    \\\mathcal T_0(\boldsymbol{a},\boldsymbol{b}) &\Leftrightarrow \rho_2(a_1,b_3)\wedge \rho_2(b_2,b_3)\wedge b_1\in\{0,1\},
    \\\mathcal T_1(\boldsymbol{a},\boldsymbol{b}) &\Leftrightarrow \rho_2(a_1,b_2)\wedge \rho_2(b_2,b_3)\wedge b_1\in\{0,1\},
    \\ \boldsymbol{a}\in U_{01} &\Leftrightarrow a_1\in\{0,1\}\wedge \rho_2(a_2,a_3).
\end{align*}
Let us define the maps $\iota$ and $\nu$ as follows:
\begin{align*}\label{equation:maph}
 \iota(x)=\begin{cases}
   (0,0,2) & \text{if } x = 0,\\
   (1,2,1) & \text{if } x = 1,\\
   (2,0,1) & \text{if } x = 2.
   \end{cases}
&&\nu(\boldsymbol{x})=\begin{cases}
   0 & \text{if } \boldsymbol{x}\in O,\\
   1 & \text{if } \boldsymbol{x}\in I,\\
   2 & \text{if } x_1 = 2.
   \end{cases}
\end{align*}
with
\begin{align*}
    O &= \Big(\{(0,b,c)\mid b,c\in\{0,1,2\}\}\setminus \{(0,2,0), (0,2,1)\}\Big) \cup \{(1,0,2), (1,1,2)\},
    \\I &= \Big(\{(1,b,c)\mid b,c\in\{0,1,2\}\}\setminus \{(1,0,2), (1,1,2)\}\Big) \cup \{(0,2,0), (0,2,1)\}.
\end{align*}

We first show that $\iota$ is a homomorphism from $\struct{M}_0'$ to $\struct{S}$. So suppose that $(a,b)\in\psi_2$ and let $(\iota(a),\iota(b)) = (\boldsymbol{a}, \boldsymbol{b})$. If $a=b=2$, then $(\iota(a),\iota(b))=((2,0,1),(2,0,1)) \in \Psi_2$. If $a\in\{0,1\}$, then $b$ must be equal to $1-a$ and $(\iota(a),\iota(b))$ is either equal to $((0,0,2),(1,2,1))$ or to $((1,2,1),(0,0,2))$. Both tuples are in $\Psi_2$. Suppose $(a,b)\in\rho_2$. If $a=2$, then $b\in\{0,1\}$. Hence $\iota(a)=(2,0,1)$ and $\iota(b)=(b_1,b_2,b_3)\in\{(0,0,2), (1,2,1)\}$. Since $(2,b_1)\in\rho_2$, we have $(\iota(a),\iota(b))\in \mathcal R_2$. By symmetry, also $(\iota(a),\iota(b))\in \mathcal  R_2$ if $a\in\{0,1\}$ and $b=2$. Suppose that $(a,b)\in\tau_0$. It follows that $b\in\{0,1\}$. If $b=1$, then $a=2$ and $(\iota(a),\iota(b))=((2,0,1),(1,2,1))\in\mathcal T_0$. Otherwise, $b=0$ and $a\in\{0,1\}$. Thus, $\iota(a)=(a_1,a_2,a_3)\in\{(0,0,2),(1,2,1)\}$ and $\iota(b)=(0,0,2)$. In both cases, it is easy to check that $(\iota(a),\iota(b))\in \mathcal T_0$.

Similarly, it is easy to check that if $(a,b)\in\tau_1$, then $(\iota(a),\iota(b))\in \mathcal T_1$. Clearly $\iota$ also preserves all unary relations, thus it is a homomorphism from $\struct{M}_0'$ to $\struct{S}$.

We next show that $\nu$ is a homomorphism from $\struct{S}$ to $\struct{M}_0'$. First suppose $(\boldsymbol{a},\boldsymbol{b})\in\Psi_2$. If $a_1=2$, then $b_1=2$ and $(\nu(a),\nu(b))=(2,2)\in\psi_2$. If $a_1\in\{0,1\}$, then $b_1= 1-a_1$. Note that $(\nu(\boldsymbol{a}),\nu(\boldsymbol{b}))\in\{(a_1,1-a_1),(1-a_1,a_1)\}$ and both tuples are in $\psi_2$. Suppose $(\boldsymbol{a},\boldsymbol{b})\in \mathcal R_2$, so $(a_1,b_1)\in\rho_2$. It follows immediately from the definition of $\nu$ that $(\nu(\boldsymbol{a}),\nu(\boldsymbol{b}))$ is in $\rho_2$. Suppose $(\boldsymbol{a},\boldsymbol{b})\in \mathcal T_0$. If $a_1=2$, then we have $\boldsymbol{b}=(b_1,2,b_3)$, with $b_1,b_3\in\{0,1\}$, i.e., $\boldsymbol{b}\in I$. Thus, $(\nu(\boldsymbol{a}),\nu(\boldsymbol{b}))=(2,1)\in\tau_0$. If $a_1\in\{0,1\}$, then $\boldsymbol{b}=(b_1,b_2,2)$, with $b_1,b_2\in\{0,1\}$, i.e., $\boldsymbol{b}\in O$. By the definition of $\nu$ it follows that $\nu(\boldsymbol{a})\in\{0,1\}$ and $\nu(\boldsymbol{b})=0$. Hence $(\nu(\boldsymbol{a}),\nu(\boldsymbol{b}))\in\tau_0$.

With a similar argument one can show that if $(a,b)\in \mathcal  T_1$, it follows that $(\nu(\boldsymbol{a}),\nu(\boldsymbol{b}))\in\tau_1$. It is also easy to see that all unary relations are preserved, thus $\nu$ is a homomorphism.

We have shown that $\mathbb S$ and $\struct{M}_0'$ are homomorphically equivalent. As $\struct{S}$ is a pp-power of $\struct{M}_0$ this finishes the proof.
\end{proof}

By the above lemmata, we have at most four $\minorequivalentto$-classes given by the clones $\mathcal I_2$, $\mathcal C_2$, $\Pol(\mathbb M_1)$, and $\Pol(\mathbb M_0)$. 

In order to show that they are all distinct, we introduce the following minor conditions.

\begin{definition} \label{def:MinorConditions2} Let us define the minor conditions $\Sigma_1$, $\Sigma_2$ as follows:
\begin{enumerate}
\item $f \models \Sigma_1$, if $f$ is a quasi majority operation that satisfies $f(x,y,z) \approx f(z,y,x)$,
\item $f \models \Sigma_2$, if $f$ is a 5-ary symmetric operation that satisfies \[f(x,x,y,y,z)\approx f(x,y,y,z,z).\]
\end{enumerate}
\end{definition}

\begin{lemma}\label{lem:separationSDM1andC2}
$\Pol(\mathbb M_1) \not\models \Sigma_1$, but $\mathcal C_2 \models \Sigma_1$.
\end{lemma}
\begin{proof}
Let $m \in \Pol(\mathbb M_1)$ be a quasi majority operation. Note that, since $\Pol(\mathbb M_1)$ is idempotent, $m$ is, in fact, a majority operation. Since $m$ preserves the congruence $\mu_2$, we find that both $m(1,2,0)$ and $m(0,2,1)$ are in the same $\mu_2$-class as $m(0,2,0) = 0$, i.e., $m(1,2,0), m(0,2,1) \in \{0,1\}$. By the preservation of $\psi_2$, it is clear that $(m(1,2,0),  m(0,2,1))\in \,\psi_2$ holds; thus, $m(1,2,0) \neq m(0,2,1)$. So $m$ cannot satisfy the identity $m(x,y,z) \approx m(z,y,x)$; hence, $\Pol(\mathbb M_1) \not\models \Sigma_1$. On the other hand, the Boolean majority operation on $\{0,1\}$ is an element of $\mathcal C_2$ that clearly witnesses $\mathcal C_2 \models \Sigma_1$.
\end{proof}

\begin{lemma}\label{lem:separationSDM0andSDM1}
$\Pol(\mathbb M_0) \not\models \Sigma_2$, but $\Pol(\mathbb M_1) \models \Sigma_2$.
\end{lemma}
\begin{proof}
    For a contradiction, let us assume that there is a 5-ary $f\in\Pol(\mathbb M_0)$ that satisfies $\Sigma_2$. Since $f$ satisfies $\Sigma_2$, $f(0,0,1,2,2) = f(0,0,2,2,1) = f(0,2,2,1,1) = f(1,1,0,2,2)$. Furthermore, since $f$ preserves $\psi_2$, we obtain
    \[\begin{pmatrix}
    f(0,0,1,2,2)
    \\f(1,1,0,2,2)
    \end{pmatrix}\in\psi_2
    \] 
   and therefore $f(0,0,1,2,2)=f(1,1,0,2,2)=2$. Since $f$ preserves $\mu_2$, this implies that $f(a,b,c,2,2)=2$, for every $a,b,c\in\{0,1\}$. Moreover,
    \[\begin{pmatrix}
    f(a,b,c,2,2)
    \\f(2,2,2,0,1)
    \end{pmatrix}=\begin{pmatrix}
    2\\d
    \end{pmatrix}\in\rho_2
    \]
    implies that $d\in\{0,1\}$. But then we would obtain
    \[\begin{pmatrix}
    f(2,2,2,0,1)
    \\f(2,2,2,1,0)
    \end{pmatrix}=\begin{pmatrix}
    d\\d
    \end{pmatrix}\in\psi_2,
    \]
which is a contradiction.

For the second part, let us define the 5-ary map 

$$f(x_1,x_2,x_3,x_4,x_5) = \begin{cases}
\mathrm{maj}_5(x_1,x_2,x_3,x_4,x_5) &\text{ if all } x_i \in \{0,1\},\\
 2 &\text { else.}
\end{cases}$$
Here $\mathrm{maj}_5$ denotes the operation on $\{0,1\}$ that returns the value that appears at least 3 times. It is not difficult to see that $f \models \Sigma_2$. As it is idempotent and preserves both $\psi_2$ and $\mu_2$, it is a polymorphism of $\mathbb M_1$.
\end{proof}

Let us summarize the results of this subsection.

\begin{theorem} \label{theorem:classesOfSDmeetMalcev}
Let $\mathcal C$ be an idempotent Mal'cev clone on $\{0,1,2\}$ that has a majority operation and preserves $\mu_2$.
\begin{enumerate}
\item If $\psi_2 \notin \Inv(\mathcal C)$, then either $\mathcal C  \minorequivalentto \mathcal C_2$ or $\mathcal C  \minorequivalentto \mathcal I_2$,
\item if $\psi_2 \in \Inv(\mathcal C)$ and $\rho_2 \notin \Inv(\mathcal C)$, then $\mathcal C  \minorequivalentto \Pol(\mathbb M_1)$,
\item if $\psi_2, \rho_2 \in \Inv(\mathcal C)$, then $\mathcal C  \minorequivalentto \Pol(\mathbb M_0)$.
\end{enumerate}
Furthermore, $\Pol(\mathbb M_0) \isaminorof \Pol(\mathbb M_1) \isaminorof \mathcal C_2 \isaminorof \mathcal I_2$ are pairwise not minor-equivalent. 
\end{theorem} 

\begin{proof}
Recall that Theorem \ref{thm:bulatovreproved} implies that 
\begin{itemize}
\item $\Pol(\mathbb H) \subseteq \mathcal C$ if and only if $\psi_2 \notin \Inv(\mathcal C)$,
\item $\psi_2 \in \Inv(\mathcal C)$ and $\rho_2 \notin \Inv(\mathcal C)$ if and only if $\Pol(\mathbb M_1'') \subseteq \mathcal C \subseteq \Pol(\mathbb M_1)$, and
\item $\psi_2, \rho_2 \in \Inv(\mathcal C)$ if and only if $ \Pol(\mathbb M_0'') \subseteq \mathcal C  \subseteq \Pol(\mathbb M_0)$.
\end{itemize}
Statement (2) and (3) then follow from Lemma~\ref{lemma:extensionm1}, and~\ref{lem:forTheLastCollapse}, respectively. Statement (1) follows from Lemma~\ref{lemma:Hppconstruction} and the fact that in $\mathfrak{P}_3$ (and even in $\Pfin$) there are no elements between $\overline{\clo{C}_2}$ and $\overline{\clo{I}_2}$~\cite{VucajZhuk, meyerSubmaximal}.

Since $\Pol(\mathbb M_0) \subseteq \Pol(\mathbb M_1) \subseteq \Pol(\mathbb H)$, we get that $\Pol(\mathbb M_0) \isaminorof \Pol(\mathbb M_1) \isaminorof \mathcal C_2$. Moreover,  $\Pol(\struct{M}_1) \notisaminorof \Pol(\struct{M}_0)$ by Lemma~\ref{lem:separationSDM0andSDM1}, $\clo{C}_2 \notisaminorof \Pol(\struct{M}_1)$ by Lemma~\ref{lem:separationSDM1andC2}, and $\clo{C}_2 \notisaminorof \clo{I}_2$ by Theorem~\ref{thm:2element}. 
\end{proof}

\section{Algebras without a majority term}\label{sec:NotSDmeet}

We are left with analyzing all idempotent Mal'cev clones that do not have a majority operation and do not preserve $\varphi$ (in particular, this also excludes the only idempotent Abelian Mal'cev clone on $\{0,1,2\}$, namely $\mathcal Z_3$, and the clone $\mathcal L_2 = \Pol(\varphi,\psi_2,T_2)$, both of which are part of the classification result for clones of self-dual operations in Theorem~\ref{thm:selfdual}).

We will show that all such clones are minor-equivalent to $\clo{Z}_2$. By the following lemma, we already see it for the minimal Mal'cev clones satisfying our assumptions.

\begin{lemma}\label{lemma:the_turning_of_the_tide}
  Let $i\in \{0,1,2\}$ and let $d_i$ be as in~\eqref{equation:thed_ifunctions}.
 Then $\clo{Z}_2 \minorequivalentto \Clo(d_i)$.
\end{lemma}
\begin{proof}
It is clear that $\Clo(d_i)\minorequivalentto\Clo(d_2)$.
 For every odd $k\geq 1$, let us define  $f_k\colon \{0,1,2\}^{k}\to \{0,1,2\}$ given for 
    all $x_1, \dots x_{k}\in \{0,1,2\}$ by 
    \[
    f_k(x_1, \dots, x_{k})=\begin{cases}
        2 & \text{ if } |\{l \colon x_{l}=2\}| \text{ is odd},\\     
        \sum_{x_j \in \{0,1\}} x_j & \text{ else},
    \end{cases}
    \]
where the latter sum is computed modulo 2. We first prove that $\Clo(d_2) = \Clo(f_k)$, for every $k \geq 3$. By definition $f_k(x,y,z,\dots,z)= d_2(x,y,z)$, and thus $\Clo(f_k) \supseteq \Clo(d_2)$. To prove the other inclusion we are going to show that $f_k$ preserves the same relations as $d_2$. By Theorem~\ref{thm:bulatovreproved}, it is enough to show that $f_k$ preserves $\{\mu_2, \rho_2, \tau_0,\tau_1, \psi_2, \psi_2' ,T_2,T_2',T^{\mu_2}, S_{01} \}$ together with all unary relations.

We start with $T_2$, so let $\boldsymbol{x}_1,\ldots,\boldsymbol{x}_{k} \in T_2$, and let $I = \{ i \in \{1,\ldots,k\} \mid \boldsymbol{x_i} = (2,2,2,2)\}$. Then, by definition,     \[
    f_k(\boldsymbol{x}_1,\ldots,\boldsymbol{x}_{k})=\begin{cases}
        (2,2,2,2) & \text{ if } |I| \text{ is odd},\\     
        \sum_{i \notin I} \boldsymbol{x}_i & \text{ else}.
    \end{cases}
    \]
So, if $|I|$ is odd, clearly $f_k(\boldsymbol{x}_1,\ldots,\boldsymbol{x}_{k}) \in T_2$. For the case that $|I|$ is even, recall that $\mathbf{u} \in T_2 \cap \{0,1\}^4$ if and only if $u_1+u_2+u_3+u_4 = 0$ modulo $2$. This also readily implies $f_k(\boldsymbol{x}_1,\ldots,\boldsymbol{x}_{k}) = \sum_{i \notin I} \boldsymbol{x}_i \in T_2$. Thus $f_k$ preserves $T_2$. It follows that also $\mu_2$ is invariant under $f_k$, as it is a coordinate kernel of $T_2$. The relation $\psi_2$ is invariant under $f_k$, since $0\mapsto 1, 1 \mapsto 0, 2\mapsto 2$ is an automorphism of $(\{0,1,2\},f_k)$. 

Next, it follows straightforward from the definition, that for all $a\neq b$ and all $x_1, \dots, x_{k} \in \{a,b\}$: 
\begin{equation} \label{eq:fkrestriction}
    f_k(x_1, \dots, x_{k})=\begin{cases}
        a & \text{ if } |\{l \colon x_{l}=a\}| \text{ is odd},\\     
        b & \text{ else},
    \end{cases}
\end{equation}
thus $f_k$ preserves every two-element set $\{a,b\}$, and therefore all unary relations. Note also that $(\{a,b\},f_k)$ is affine, and has the automorphism $a \mapsto b, b \mapsto a$. For $\{a,b\} = \{0,1\}$, this implies that $f_k$ preserves $T_2'$ and $\psi_2'$. The quotient algebra $(\{0,1,2\},f_k)/\mu_2$ also satisfies the equation in \eqref{eq:fkrestriction}, so in the same way we see that $f_k$ preserves $T^{\mu_2}$ and $\rho_2$. The two isomorphisms between $(\{0,1,2\},f_k)/\mu_2$ and $(\{0,1\},f_k|_{\{0,1\}})$ correspond to the relations $\tau_0$ and $\tau_1$ being invariant under $f_k$. Last, $S_{01}$ is pp-definable from the above relations by $\exists y \, \tau_0(x_3,y) \land T_2(x_1,x_2,y,0)$. Hence $\Clo(f_k) = \Clo(d_2)$, for every odd $k\geq 3$.

In order to prove that $\Clo(d_2) \minorequivalentto \clo{Z}_2$, first recall that the restriction of $(\{0,1,2\},d_2)$ to $\{0,1\}$ is $\mathbf Z_2 =(\{0,1\},x+y+z)$. Thus the restriction map $f \mapsto f|_{\{0,1\}}$ is a minion homomorphism from $\Clo(d_2)$ to $\clo{Z}_2$ (in fact, even a clone homomorphism).

In order to find a minion homomorphism from $\clo{Z}_2$ to $\Clo(d_2)$, let us first define, for every $l \in \N$ and $I \subseteq [l]$ of odd cardinality the $l$-ary operation $f_I(x_1,\ldots,x_l) = f_{|I|}(x_{i_1},\ldots,x_{i_{|I|}})$. Note that this does not depend on the enumeration of $I$, since $f_k$ is symmetric.

The minor identities $f_{k+2}(y,y,x_1,\ldots,x_k) \approx f_{k}(x_1,\ldots,x_{k})$ hold, for every $k\geq 3$. Thus, whenever we consider a minor $(f_k)_\alpha(x_1,\ldots,x_l) = f_k(x_{\alpha(1)},\ldots x_{\alpha(k)})$, for a map $\alpha \colon [k] \to [l]$, the above identities imply that $(f_k)_\alpha = f_I$, where $i \in I$ iff $\alpha^{-1}(i)$ is of odd cardinality. Note that $|I|$ must be odd, since $k$ is odd. Thus, the closure of $\{f_1,f_3,\ldots\}$ under forming minors consists exactly of the maps $f_I$.

Next recall that every $k$-ary function in $\clo{Z}_2$ is of the form $b_I(\boldsymbol{x}) = \sum_{i\in I} x_i$, where $I \subseteq [k]$ is of odd size. It is then follows straightforwardly from the above discussion, that the map $\xi\colon \clo{Z}_2 \to \Clo(d_2)$ given by $\xi(b_I) = f_I$ is a minion homomorphism.

We remark that, more generally it is know that $\clo{Z}_2$ has a minion homomorphism to $\Pol(\mathbb A)$ for a finite structure $\mathbb A$, if and only if $\Pol(\mathbb A)$ contains symmetric terms $f_k$ of all odd arities satisfying $f_{k}(x_1,x_2,\ldots,x_k,y,y) \approx f_k(x_1,x_2,\ldots,x_k,z,z)$, for every $k$. This folklore result can be proven by a slight modification of the argument in \cite[Theorem 7.19.]{jakubPCSP}.

\end{proof}

By the next lemma, any other Mal'cev clone satisfying our assumptions must contain a minimal Mal'cev clone of the form $\Clo(d_i)$.

\begin{lemma}\label{lemma:algebras_not_sdmeet_di}
Let $\clo{C}$ be an idempotent Mal'cev clone on $\{0,1,2\}$ such that $\clo{C}$ has no majority operation and does not preserve $\varphi$. There exists $i\in \{0,1,2\}$ such that $d_i\in \clo{C}$.
\end{lemma}
\begin{proof}
First, assume that $\clo{C}$ preserves a nontrivial equivalence relation $\mu$; without loss of generality, $\mu = \mu_2$. Then, by Theorem~\ref{thm:bulatovreproved}, there exists a relational basis of $\clo{C}$ that consists of a subset of 
$\{\mu_2, \rho_2, \tau_0,\tau_1, \psi_2, \psi_2' ,T_2,T_2',T^{\mu_2}, S_{01} \}$ together with unary relations. Since $d_2$ preserves all such relations, we get $d_2\in \clo{C}$.

If $\clo{C}$ preserves no nontrivial equivalence relation, then (up to renaming elements), by Theorem~\ref{thm:bulatovreproved}, there exists a relational basis consisting of a subset of 
$\{\varphi, \psi_2, \psi_0', \psi_1', \psi_2', \varphi_0', \ldots, \varphi_5',$ $ T,T_1',T_2',T_3'\}$ and unary relations. 
Since $\clo{C}$ does not preserve $\varphi$, it also cannot preserve $T$, as $(x,y) \in \varphi \Leftrightarrow T(x,1,y,0)$. Thus, 
a relational basis exists that is a subset of  $\{\psi_2, \psi_0', \psi_1', \psi_2', \varphi_0', \ldots, \varphi_5',T_1',T_2',T_3'\}$ and unary relations. All such relations are preserved by $d_2$; thus, $d_2 \in \clo{C}$.
\end{proof}

Lemma~\ref{lemma:algebras_not_sdmeet_di} already implies that all clones satisfying our assumption lie above~$\clo{Z}_2$. To show equivalence, we need the following lemma.

\begin{lemma}\label{lemma:algebras_not_sdmeet_subalgebras}
    Let $\algA$ be an idempotent Mal'cev algebra on $\{0,1,2\}$ such that $\varphi \notin \Inv(\algA)$, and $\alg{A}$ does not have a majority term. Then $\alg{A}$ has a two-element Abelian subalgebra or a two-element Abelian proper quotient. 
\end{lemma}
\begin{proof}
By Theorem~\ref{thm:majority}, there exists an algebra $\algB \in \HS(\algA)$ such that $\algB$ has an Abelian monolith. If $\algB$ is a 2-element algebra, then its monolith must be the full relation, so $\algB$ is Abelian and we are done. Thus, let us assume that $\algB = \algA$ has an Abelian monolith. Note that $\algA$ cannot be simple, as then $\Clo(\algA) = \mathcal Z_3$, which preserves $\varphi$. Thus, $\algA$ has a proper Abelian monolith $\mu$. Since $\algA$ is idempotent, every equivalence class of $\mu$ is also a subalgebra; this gives us an Abelian 2-element subalgebra.
\end{proof}

We conclude:

\begin{theorem}\label{thm:notmajority}
Let $\clo{C}$ be an idempotent Mal'cev clone on $\{0,1,2\}$ such that $\clo{C}$ has no majority operation and $\varphi \notin \Inv(\clo{C})$. Then $\clo{C}\minorequivalentto \clo{Z}_2$. 
\end{theorem}
\begin{proof}
Let $\algA$ be an algebra with $\Clo(\algA) = \clo{C}$. By Lemma~\ref{lemma:algebras_not_sdmeet_subalgebras}, $\alg{A}$ has either a two-element Abelian subalgebra or a two-element Abelian proper quotient $\algB$. Since $\algB$ is idempotent, $\Clo(\algB)$ needs to be equal to $\clo{Z}_2$, the clone of idempotent affine operations modulo $2$. The map $t^{\algA} \mapsto t^{\algB}$ naturally yields a clone homomorphism from $\Clo(\algA) = \clo{C}$ to $\Clo(\algB) = \clo{Z}_2$; thus, $\clo{C} \isaminorof \clo{Z}_2$. 

On the other hand, Lemma~\ref{lemma:algebras_not_sdmeet_di} together with Lemma~\ref{lemma:the_turning_of_the_tide} implies that $\mathcal Z_2 \minorequivalentto\Clo(d_i) \subseteq \mathcal C$; thus, $\clo{C}\minorequivalentto \clo{Z}_2$.
\end{proof}

\section{Main result}\label{sec:Main}
In this section, we prove the characterization of Mal'cev clones on $\{0,1,2\}$ modulo minor-equivalence by gathering together the results from the previous sections. Recall that in Section~\ref{sec:selfdual}, we already identified 8 elements; in Section~\ref{sec:SDmeet}, we further introduced the clones $$ \clo{M}_0 = \Pol(\struct{M}_0) \text{ and } \clo{M}_1 = \Pol(\struct{M}_1).$$
We claim that these 10 clones form representatives of all equivalence classes of quasi Mal'cev clones in $\Pthree$. In the table in Figure~\ref{figure:definitionclones}, we recall their definitions. As the first four clones were defined on sets of size $2$ or smaller, we also provide a $\minorequivalentto$-representative on $\{0,1,2\}$ for them in the third column. In the case of $\clo{C}_2$ and $\clo{Z}_2$, these representatives follow from Lemma~\ref{lem:C2Collapse} and Theorem~\ref{thm:notmajority}, respectively. A complete list of all representatives for each $\minorequivalentto$-class can be found in the database provided in the Appendix of this article.

\begin{figure}[h]
\begin{tabular}{|l|l|l|}
\hline
Clone & Definition & $\minorequivalentto$ equivalent to \\
\hline
$\clo{T}$ & $\Pol(\{0\};\emptyset)$ & $\Pol(\{0,1,2\};\emptyset)$\\
$\clo{I}_2$ & $\Pol(\{0,1\};\{0\},\{1\})$, & $\Pol(\{0\},\{1\},\{2\})$\\
$\clo{C}_2$ & $\Pol(\mathbb C_2)$ & $\Pol(\psi'_2,\{0\},\{1\},\{2\})$\\
$\clo{Z}_2$ & $\Pol(\struct{Z}_2)$ & $\Pol(\psi'_2,T_2')$\\
$\clo{C}_3$ & $\Pol(\mathbb C_3)$ & \\
$\clo{M}_1$ & $\Pol(\struct{M}_1)$ & \\
$\clo{M}_0$ & $\Pol(\struct{M}_0)$ &\\
$\clo{D}$ & $\Pol(\{0,1,2\};\varphi,\psi'_2)$ & \\
$\clo{L}_2$ &  $\Pol(\{0,1,2\};\varphi,\psi_2,T_2)$ & \\
$\clo{Z}_3$ & $\Pol(\{0,1,2\};T, \{0\},\{1\},\{2\})$  &\\
\hline
\end{tabular}
    \caption{Relational description of representatives of all equivalence classes of Mal'cev clones in $\Pthree$}
    \label{figure:definitionclones}
\end{figure}

In order to get a full description of the order induced by $\isaminorof$ on these 10 clones, we still need the following observations:

\begin{lemma}\label{lem:Separations}
    The following hold:
    \begin{enumerate}
        \item $\clo{L}_2 \isaminorof \clo{Z}_2  \isaminorof \clo{M}_0$ and $\clo{L}_2 \isaminorof \clo{D}$, \label{itm:inclusion}
        \item $\clo{L}_2\notisaminorof\clo{Z}_3$, by the existence of a (quasi) minority operation, \label{itm:Sep1}
        \item $\clo{Z}_3\notisaminorof\clo{L}_2$ and $\clo{C}_3\notisaminorof\clo{Z}_2$, by the existence of a 2-cyclic operation,\label{itm:Sep2}
        \item $\clo{Z}_2\notisaminorof\clo{C}_3$, by the existence of a 3-cyclic term,\label{itm:Sep3}
        \item $\clo{M}_0\notisaminorof\clo{Z}_2$ and $\clo{D}\notisaminorof\clo{L}_2$, by the existence of a (quasi) majority operation.\label{itm:Sep4}
    \end{enumerate}
\end{lemma}
\begin{proof}
    \eqref{itm:inclusion} All of the statements follow from our results with comparisons up to term inclusion. $\clo{L}_2 = \Pol(\varphi,\psi_2,T_2) \subseteq \Pol(\psi_2,T_2) \subseteq \Pol(\mathbb M_0) = \clo{M}_0$. The second clone is equivalent to $\clo{Z}_2$ by Theorem~\ref{thm:notmajority}. Similarly, the inclusion $\clo{L}_2 \subseteq \clo{D}$ clearly implies $\clo{L}_2 \isaminorof \clo{D}$.
    
    \eqref{itm:Sep1} The operation $g\in\clo{L}_2$ (see \eqref{equation:minority_plu_porj}) is a minority operation. If $\clo{Z}_3$ had a minority operation $f$, such operation would not preserve the relation $T$ since
    \[\begin{pmatrix}
        f(0,0,1)
        \\f(0,1,0)
        \\f(0,0,1)
        \\f(0,2,2)
    \end{pmatrix}=\begin{pmatrix}
       1 \\1\\1\\0
    \end{pmatrix}\notin T.\]
    As both $\clo{L}_2$ and $\clo{Z}_3$ are idempotent, we get $\clo{L}_2\notisaminorof\clo{Z}_3$, since $\clo{L}_2$ has a quasi minority, but $\clo{Z}_3$ does not (note that this was already pointed out in \cite{BodirskyVucajZhuk}, see Figure~\ref{fig:selfdual}).
    
    \eqref{itm:Sep2} The operation $2(x+y)\mod3$ is in the clone $\clo{Z}_3$ and is a 2-cyclic operation; however, $\clo{L}_2$ does not have a 2-cyclic operation, as this contradicts the preservation of $\psi_2'$. Hence $\clo{Z}_3\notisaminorof\clo{L}_2$.

    Note that $2(x+y)\bmod3$ is also in $\clo{C}_3$, but $\clo{Z}_2$ does not have a 2-cyclic operation, as it preserves inequality. Hence $\clo{C}_3 \notisaminorof \clo{Z}_2$.

      \eqref{itm:Sep3} $\clo{Z}_2$ clearly has a 3-cyclic term, e.g., the Boolean minority operation $\minority$. However, the preservation of $\varphi$ implies that $\clo{C}_3$ does not have a 3-cyclic term, thus $\clo{Z}_2\notisaminorof\clo{C}_3$. 

    \eqref{itm:Sep4} For this, simply recall that $\clo{Z}_2$ and $\clo{L}_2$ do not have a (quasi) majority operation, but $\clo{M}_0$ and $\clo{D}$ do.
\end{proof}

\begin{figure}[ht]
    
    \centering\includegraphics{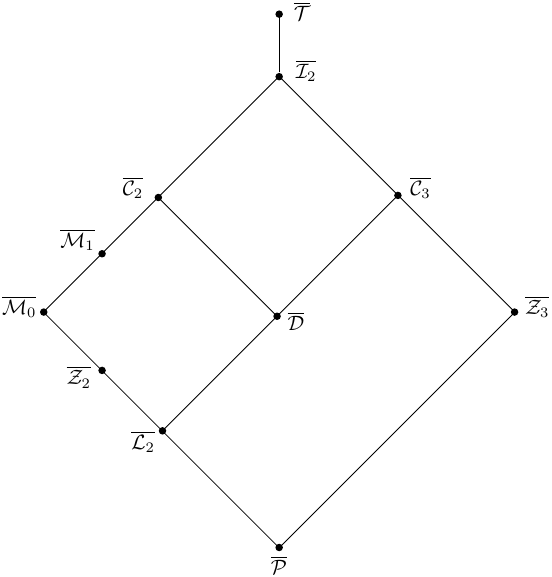}
    \caption{A Hasse diagram of the 10 element of $\Pthree$ that have a quasi Mal'cev operation (and, additionally, the bottom element $\overline{\mathcal P})$ of $\Pthree$}
    \label{figure:main_figure}

\end{figure}

We are now ready to prove our main result:

\begin{theorem}\label{teor:each_malcev_clone_in_one_class}
    Let $\clo{C}$ be a clone with a quasi Mal'cev operation on a three-element set. Then there exists 
    \[
    \clo{A}\in 
    \{\clo{Z}_3,
    \clo{L}_2,
    \clo{Z}_2,
    \clo{M}_0,
    \clo{M}_1,
    \clo{C}_2,
    \clo{D},
    \clo{C}_3,
    \clo{I}_2,
    \clo{T}\}
    \]
    such that $\clo{C}\minorequivalentto\clo{A}$. The order $\isaminorof$ on the $\minorequivalentto$-classes is depicted in the Hasse diagram in Figure~\ref{figure:main_figure}. A list of minor conditions separating different classes is given in Figure~\ref{fig:Separations}.
\end{theorem}


\begin{figure}[ht]

\centering\begin{tabular}{r|llllllllll}
& $\clo{Z}_3 \not \models$ & $\clo{L}_2 \not \models$ & $\clo{Z}_2 \not \models$ & $\clo{M}_0 \not \models$ & $\clo{M}_1 \not \models$ & $\clo{D} \not \models$ & $\clo{C}_2 \not \models$ & $\clo{C}_3 \not \models$ & $\clo{I}_2 \not \models$  \\
\hline
$\clo{Z}_3 \models $ & & 2-cyc & 2-cyc & 2-cyc & 2-cyc & 2-cyc & 2-cyc &  & & \\
$\clo{L}_2 \models $ & min & &  &  &  &  &  &  &  & \\
$\clo{Z}_2 \models $ & min & 3-cyc & &  &  & 3-cyc &  & 3-cyc &  & \\
$\clo{M}_0 \models $ & min & 3-cyc & maj &  &  & 3-cyc &  & 3-cyc &  & \\
$\clo{M}_1 \models $ & min & 3-cyc & maj & $\Sigma_2$ &  & 3-cyc &  & 3-cyc &  & \\
$\clo{D} \models $ & min & maj & maj & $\Sigma_1$ & $\Sigma_1$ &  &  &  &  & \\
$\clo{C}_2 \models $ & min & 3-cyc & maj &  $\Sigma_1$ &  $\Sigma_1$ & 3-cyc &  & 3-cyc &  & \\ 
$\clo{C}_3 \models $ & min & maj & maj & 2-cyc & 2-cyc & 2-cyc & 2-cyc &  &  & \\
$\clo{I}_2 \models $ & min & 3-cyc & maj & 2-cyc & 2-cyc & 2-cyc & 2-cyc & 3-cyc &  & \\
$\clo{T} \models $ & min & 3-cyc & maj & 2-cyc & 2-cyc & 2-cyc & 2-cyc & 3-cyc & $f(x) \approx f(y)$ & 
\end{tabular}

\caption{Minor conditions witnessing $\notisaminorof$. Here \emph{min}, \emph{maj}, and \emph{n-cyc} indicate the existence of a \emph{quasi minority}, \emph{quasi majority}, and \emph{n-cyclic} operation, respectively (see Definition~\ref{def:MinorConditions}). $\Sigma_1$ and $\Sigma_2$ were defined in Definition~\ref{def:MinorConditions2}.}
\label{fig:Separations}
\end{figure}

\begin{proof}[Proof of Theorem \ref{teor:each_malcev_clone_in_one_class}]
Let us fix some relational structure $\mathbb B$, such that $\mathcal C = \Pol(\mathbb B)$. By $\mathbb A$, let us denote the core of $\mathbb B$, extended by all singleton unary relations. Then, by Proposition~\ref{prop:idempotent}, we know that $\mathcal C \minorequivalentto \Pol(\mathbb A)$. Note that $\Pol(\mathbb A)$ contains a Mal'cev operation, as it is idempotent.

We distinguish several cases: 

\begin{description}
\item[Case~1] $|A| \leq 2$. Then, by Theorem~\ref{thm:2element}, we get that $\clo{C}$ is equivalent to an element from $\{\clo{Z}_2,\clo{C}_2,\clo{I}_2,\clo{T}\}$.
\item[Case~2] $|A| = 3$. By Theorem~\ref{thm:bulatovreproved}, without loss of generality, assume that $\mathbb A$ only contains the, at most 4-ary, standard relations (see \eqref{eq:standardrel} and \eqref{eq:standardrel2}), and all singleton unary relations. Without loss of generality, let us further assume that $\mathcal C = \Pol(\mathbb A)$.
\begin{description}
\item[Case~2a] $\mathcal C$ preserves $\varphi$, i.e., it consists of self-dual operations. In this case, by Theorem~\ref{thm:selfdual}, $\clo{C}$ is equivalent to an element of $\{\clo{C}_3, \clo{D}, \clo{L}_2, \clo{Z}_3\}$. 

\item[Case~2b] $\mathcal C$ does not preserve $\varphi$, but contains a majority operation. By Lemma~\ref{lemma:simple} and  Theorem~\ref{theorem:classesOfSDmeetMalcev}, $\clo{C}$ is equivalent to an element of $\{\clo{I}_2,\clo{C}_2,\clo{M}_1,\clo{M}_0\}$.

\item[Case~2c] $\mathcal C$ does not preserve $\varphi$ and does not contain a majority operation. Theorem~\ref{thm:notmajority} implies that $\clo{C} \minorequivalentto \clo{Z}_2$.
\end{description}
\end{description}

This finishes the proof that $\clo{C}$ is minor-equivalent to an element of $\{\clo{Z}_3, \clo{L}_2, \clo{Z}_2, \clo{M}_0,\clo{M}_1, \clo{C}_2,$ $ \clo{D},\clo{C}_3, \clo{I}_2,\clo{T}\}$.

The fact that the order depicted in Figure~\ref{figure:main_figure} is equal to $\isaminorof$ follows directly from the existing classifications in Theorem~\ref{thm:2element}, Theorem~\ref{thm:selfdual}, and our results in Theorem~\ref{theorem:classesOfSDmeetMalcev} and Lemma~\ref{lem:Separations}. In Figure~\ref{fig:Separations}, we give a list of minor conditions that can be used to see that $\mathcal A\notisaminorof \mathcal B$, for all pairs $(\mathcal A, \mathcal B)$ in our classification that are not $\isaminorof$-related. These conditions  are again based on the classifications in Theorem~\ref{thm:2element}, Theorem~\ref{thm:selfdual}, as well as Lemma~\ref{lem:separationSDM0andSDM1}, Lemma~\ref{lem:separationSDM1andC2}, and Lemma~\ref{lem:Separations} (2)-(4). Moreover, the $\minorequivalentto$-class of the projection clone $\clo{P}$ clearly lies strictly below all 10 classes of Mal'cev clones, and is the only element of the order in Figure~\ref{figure:main_figure} that does not have a quasi Mal'cev operation. This finishes our proof.
\end{proof}

Note that we included the bottom element of $\Pthree$, the equivalence class of the projection clone $\clo{P}$ in Figure~\ref{figure:main_figure}, although it is not a Mal’cev clone. Nevertheless, we decided to do so to emphasize that $\overline{\clo{L}_2}$ and $\overline{\clo{Z}_3}$ are indeed atoms of $\Pthree$, see~\cite{MinTaylorOver3}.

\section{Concluding Remarks and Future Directions} \label{sec:future}

We conclude this article with some open problems. As already mentioned in the introduction, in order to obtain a complete description of $\Pthree$, it remains to describe the downset of $\overline{\clo{B}_2}$:

\begin{problem} \label{prob:Prob1}
Describe all the elements of $\Pthree$ that are smaller than $\overline{\clo{B}_2}$, with respect to~$\isaminorof$.
\end{problem}

A major challenge in solving problem~\ref{prob:Prob1} is that there is no explicit description of the downset of $\clo{B}_2$ in $\mathfrak L_3$ (in contrast to the ideal of self-dual clones described in \cite{Zhuk15}). To date, the cardinality of $\Pthree$ remains unknown, and the same holds for $\Pfin$. From \cite[Theorem 2.13]{VucajZhuk} and the main result from \cite{AichingerMayrMcKenzie} it follows that if $\Pthree$ (or $\Pfin$) does indeed have continuum many elements, then these must necessarily belong to the ideal generated by $\overline{\clo{B}_2}$.

Another promising direction would be to extend the work done in this article to larger finite domains or, more ambitiously, to all finite sets. More precisely, one could aim to classify Mal'cev clones on arbitrary finite sets up to minion homomorphisms. 

\begin{problem}\label{prob:Prob2}
Describe all Mal'cev clones over some finite set up to minion homomorphisms.
\end{problem}

To our current knowledge, we do not know of any constructive way to obtain a relational basis for a given Mal'cev algebra; thus, completely solving Problem~\ref{prob:Prob2} might be a too ambitious project. However, in the following, we outline a few natural starting points.

As we discussed in Theorem~\ref{thm:majority}, a Pixley clone, i.e., a clone that has both a Mal'cev and a majority operation, always has a binary relational basis. This leads to the question:

\begin{problem} \label{prob:Prob3}
Describe all Pixley clones over some finite set up to minion homomorphisms.
\end{problem}

Significant progress on Problem~\ref{prob:Prob3} was already made in \cite{meyerSubmaximal}, where Meyer and Starke classify all the submaximal elements of $\Pfin$. From their analysis, it follows that, besides $\overline{\mathcal B_2}$, all such submaximal elements are the ($\minorequivalentto$-classes of) polymorphism clones of the Cayley graphs of finite simple groups $\alg{G}$. In particular, all such clones are Pixley clones.

In a recent publication, Bodirsky and Moorhead \cite{BM-conservativemaltsev} also discussed conservative Mal'cev clones. Such clones are more feasible to work with, as they always have an most 4-ary relational basis. Although \cite{BM-conservativemaltsev} does not give a full classification of conservative Mal'cev clones up to $\isaminorof$, it describes minimal elements and, based on this analysis, shows that all corresponding CSPs are in the complexity class $\bigoplus L$.

A natural class, which subsumes all 3-element Mal'cev clones, Pixley clones, and conservative Mal'cev clones, are Mal'cev clones, whose associated algebras generate residually finite varieties. Thus, we ask:

\begin{problem} \label{prob:Prob4}
Describe $\isaminorof$ for all clones of finite Mal'cev algebras that generate a residually finite variety. What are the minimal elements, and what is the computational complexity of the corresponding CSPs?
\end{problem}

In Theorem \cite[Theorem 4.3]{KS-parallelogramterms}, Kearnes and Szendrei gave finite upper bounds on the relational basis of such algebras. However, in contrast to arity 4, which suffices for all mentioned examples, their upper bound is exponential in the size of the domain of the clone. Thus, any progress on Problem~\ref{prob:Prob4} would probably require a more in-depth analysis of the relational bases in this case.

\appendix

\section{Appendix}\label{sec:appendix}

Under the following GitHub link\footnote{\href{https://github.com/StefanoF66/Mal-cev-clones-over-0-1-2-up-to-minor-equivalence.git}{https://github.com/StefanoF66/Mal-cev-clones-over-0-1-2-up-to-minor-equivalence.git}}, we present a repository containing databases related to classification presented in the paper. The data is organized into three directories, each providing different aspects of the classification. 

\subsection*{1. Standard Relations}  
The \texttt{Standard\_relations} directory contains a JSON file that lists all relations used in this paper and by A. Bulatov in~\cite{BulatovMaltsev} to characterize Mal'cev clones on a 3-element set using the notation of the current paper.

\subsection*{2. Functional Generators}  
The \texttt{Functional\_generators} directory includes data on functional generators of polymorphism clones, structured as follows:

\begin{itemize}
    \item \texttt{algebras\_with\_a\_monolith.csv} and \texttt{simple\_algebras.csv} list a set of functions whose power set characterizes all Mal'cev clones on a 3-element set. Each table provides:
    \begin{itemize}
        \item The functional generators.
        \item The standard relations they preserve.
    \end{itemize}
    To extract the functional generators of \(\Pol(\mathbb R)\) for a given structure \( \mathbb R = (A; R_1, \dots, R_n) \), query the dataset and select all rows where the value in the columns corresponding to \( R_i \) (for \( i = 1, \dots, n \)) is 1.

    \item \texttt{unary\_functions.csv}, \texttt{binary\_functions.csv}, and \texttt{3-ary\_and\_4-ary} \texttt{\_functions.csv} are CSV files with the multiplication tables of all functional generators, categorized by arity. The last file includes the only $4$-ary functional generator, which is represented as a ternary function fixing its first component. Specifically, the multiplication table of $f232\_i$ corresponds to the values of $f232$ evaluated on tuples $(i, x_1, x_2, x_3)$ for $x_1, x_2, x_3 \in \{0,1,2\}$.

\end{itemize}

\subsection*{3. Structure Databases}  
The \texttt{Structures\_databases} directory contains the relational structure corresponding to all Mal'cev clones on a 3-element set, categorized by the congruence lattices of the associated algebra. In cases where a nontrivial \textbf{monolith} exists, it is always assumed to be \( \{\{0,1\}, \{2\}\} \) and is not explicitly listed. In the simple case, only the structures with non-Abelian algebras associated with the polymorphism clones are listed since those with an Abelian one can be trivially detected.

\begin{itemize}
    \item For each congruence lattice type, two files are provided:
    \begin{enumerate}
        \item A list of structures with:
        \begin{itemize}
            \item standard relations included,
            \item unary polymorphisms,
            \item core cardinality,
            \item minor-equivalence classes of the cores (when the structure is not already a core),
            \item minor conditions satisfied (following the characterization of Figure~\ref{fig:Separations}),
            \item minor-equivalence classification.
        \end{itemize}
        \item The first list restricted to 3-element cores. To further restrict this list to the idempotent cores (the ones that are most relevant for our paper), further filter the database by selecting the rows with all 1-element subalgebras. 
    \end{enumerate}
\end{itemize}

\section*{Acknowledgments}
This paper was supported by Charles University under grant number PRIMUS/24/SCI/008. Michael Kompatscher was additionally funded by the Czech Science Foundation (GA\v{C}R grant no. 25-16324S) and the Charles University Research Center program No. UNCE/24/SCI/022.

\section*{ORCID}
\noindent Stefano Fioravanti - \url{https://orcid.org/0000-0001-6918-1805}

\noindent Michael Kompatscher - \url{https://orcid.org/0000-0002-0163-6604}

\noindent Bernardo Rossi - \url{https://orcid.org/0000-0002-0404-6634}

\noindent Albert Vucaj - \url{https://orcid.org/0000-0001-8394-9745}

\bibliographystyle{plain}
\bibliography{main.bib}

@preamble{"\def\cprime{$'$} "}

@Article{Post,
author = {Emil L. Post}, 
title = {The two-valued iterative systems of mathematical logic}, 
journal = {Annals of Mathematics Studies}, 
volume = {5}, 
publisher = {Princeton University Press},
address = {Princeton},
year = {1941}}

@Article{3elem,
author = {Yu. I. Yanov and A. A. Muchnik},
title = {On the existence of $k$-valued closed classes that have no bases}, 
journal = {Dokl. Akad. Nauk SSSR},
volume = {127},
year = {1959}, 
pages = {44-46}}

@article{JabMaximal,
author = {S. V. Jablonskij},
title = {On functional completeness in the three-valued calculus}, 
journal = {Dokl. Akad. Nauk SSSR}, 
volume = {95(2)},
pages = {1153–1155}, 
year = {1954},
note = {(in Russian)}
}

@article{Csakany84,
  author    = {B\'ela Cs{\'a}k{\'a}ny},
  title     = {All minimal clones on the three-element set},
  journal   = {Acta Cybern.},
  volume    = {6},
  year      = {1984},
  pages     = {227-238},
  bibsource = {DBLP, http://dblp.uni-trier.de}
}

@Article{Rosenberg,
author = {Ivo G. Rosenberg},
title = {Minimal clones {I}: the five types}, 
journal = {Lectures in Universal Algebra 
(Proc. Conf. Szeged, 1983), Colloq. Math. Soc. J. Bolyai},
year = {1986},
volume = {43},
pages = {405--427}
}

@Article{DemetrHannak,
author = {J. Demetrovics and L. Hannak},
title = {The cardinality of selfdual closed classes in $k$-valued logics}, 
journal = {MTA SzTAKI K\"ozlemenyek}, 
volume = {23}, 
pages = {8-17}, 
year = {1979}}

@article{RosenbergMaximal,
author = {Ivo~G. Rosenberg},
title = {La structure des fonctions de plusieurs variables sur un ensemble
  fini},
journal = {C. R. Acad. Sci., Paris},
volume = {260},
pages = {3817-3819},
year = {1965}
}

@article{Marchenkov,
author = {Sergey S. Marchenkov},
title = {On closed classes of self-dual functions of many-valued logic},
journal = {Probl. Kibernetiki},
pages = {261--266},
volume = {40},
year = {1983},
note = {(in Russian)}
}

@article{Zhuk15,
  author    = {Dmitriy Zhuk},
  title     = {The Lattice of All Clones of Self-Dual Functions in Three-Valued Logic},
  journal   = {Multiple-Valued Logic and Soft Computing},
  volume    = {24},
  number    = {1-4},
  pages     = {251--316},
  year      = {2015}
}

@article {BulatovImportante,
    AUTHOR = {Bulatov, Andrei A.},
     TITLE = {Three-element {M}al'tsev algebras},
    JOURNAL = {Acta Scientiarum Mathematicarum},
    VOLUME = {71},
      YEAR = {2005},
    NUMBER = {3-4},
     PAGES = {469--500},
      ISSN = {0001-6969,2064-8316},
   MRCLASS = {08A40 (08B05)},
  MRNUMBER = {2206592},
MRREVIEWER = {Keith\ A.\ Kearnes},
}

@article{KS-parallelogramterms,
  title={Clones of algebras with parallelogram terms},
  author={Kearnes, Keith A. and Szendrei, {\'A}gnes},
  journal={International Journal of Algebra and Computation},
  volume={22},
  number={01},
  pages={1250005},
  year={2012},
  publisher={World Scientific},
  doi={10.1142/S0218196711006716}
}

@book{garciaTaylor,
  title={The Lattice of Interpretability Types of Varieties},
  author={Garcia, O.C. and Taylor, W.},
  isbn={9780821839591},
  series={Memoirs of the American Mathematical Society Series},
  url={https://books.google.de/books?id=2ESoPwAACAAJ},
  year={2005},
  publisher={American Mathematical Society}
  }

@Article{wonderland,
author={Libor Barto and Jakub Opr\v{s}al and Michael Pinsker},
title={The wonderland of reflections},
journal = {Israel Journal of Mathematics},
volume=223,
number=1,
year=2018,
pages={363-398}
}

@inproceedings{BulatovFVConjecture,
  author    = {Andrei A. Bulatov},
  title     = {A Dichotomy Theorem for Nonuniform {CSP}s},
  booktitle = {58th {IEEE} Annual Symposium on Foundations of Computer Science, {FOCS}
               2017, Berkeley, CA, USA, October 15-17, 2017},
  pages     = {319--330},
  year      = {2017}
}

@inproceedings{ZhukFVConjecture,
  author    = {Dmitriy Zhuk},
  title     = {A Proof of {CSP} Dichotomy Conjecture},
  booktitle = {58th {IEEE} Annual Symposium on Foundations of Computer Science, {FOCS}
               2017, Berkeley, CA, USA, October 15-17, 2017},
  pages     = {331--342},
  year      = {2017}
}

@article{ZhukDichotomy, 
author = {Zhuk, Dmitriy}, 
title = {A Proof of the {CSP} Dichotomy Conjecture}, 
year = {2020}, 
issue_date = {October 2020}, 
journal = {J. ACM}, 
address = {New York, NY, USA}, 
volume = {67}, 
number = {5}, 
issn = {0004-5411}, 
url = {https://doi.org/10.1145/3402029}, 
doi = {10.1145/3402029}}

@ARTICLE{jakubPCSP,
       author = {{Barto}, Libor and {Bul{\'\i}n}, Jakub and {Krokhin}, Andrei and
         {Opr{\v{s}}al}, Jakub},
        title = "{Algebraic approach to {P}romise {C}onstraint {S}atisfaction}",
  journal={Journal of the ACM (JACM)},
  volume={68},
  number={4},
  pages={1--66},
  year={2021},
  publisher={ACM New York, NY}
}

@phdthesis {FlorianThesis,
    author = {Florian Starke}, 
     title = {Digraphs modulo primitive positive constructability}, 
    school = {Technische Universität Dresden}, 
      year = {2024},
}

@phdthesis {AlbertThesis,
    author = {Albert Vucaj}, 
     title = {Clones over Finite Sets and Minor Conditions}, 
    school = {Technische Universität Dresden}, 
      year = {2023},
}

@article{cyclesBodirskyStrakeVucaj,
author = {Bodirsky, Manuel and Starke, Florian and Vucaj, Albert},
title = {Smooth digraphs modulo primitive positive constructability and cyclic loop conditions},
journal = {International Journal of Algebra and Computation},
volume = {31},
number = {05},
pages = {929-967},
year = {2021},
doi = {10.1142/S0218196721500442},

URL = { 
        https://doi.org/10.1142/S0218196721500442
    },
}

@article{VucajZhuk,
    AUTHOR = {Vucaj, Albert and Zhuk, Dmitriy},
     TITLE = {Submaximal clones over a three-element set up to
              minor-equivalence},
   JOURNAL = {Algebra Universalis},
    VOLUME = {85},
      YEAR = {2024},
    NUMBER = {2},
     PAGES = {Paper No. 22, 31},
      ISSN = {0002-5240,1420-8911},
   MRCLASS = {03B50 (08A70 08B05)},
  MRNUMBER = {4719392},
       DOI = {10.1007/s00012-024-00852-w},
       URL = {https://doi.org/10.1007/s00012-024-00852-w},
}

@unpublished{meyerSubmaximal,
      title={Finite Simple Groups in the Primitive Positive Constructability Poset}, 
      author={Sebastian Meyer and Florian Starke},
      year={2024},
      eprint={2409.06487},
      archivePrefix={arXiv},
      url={https://arxiv.org/abs/2409.06487}, 
      note={arXiv preprint arXiv:2409.06487}
}

@article{BodirskyStarke,
      title={Maximal Digraphs With Respect to Primitive Positive Constructibility}, 
      author={Manuel Bodirsky and Florian Starke},
      journal={ Combinatorica},
      year={2022}, DOI={https://doi.org/10.1007/s00493-022-4918-1}
}

@ARTICLE{albert,
  author = {Bodirsky, Manuel and Vucaj, Albert}, 
  title = {{Two-element structures modulo primitive positive constructability}}, 
  issn = {0002-5240}, 
  doi = {10.1007/s00012-020-0647-8}, 
  pages = {20}, 
  number = {2}, 
  volume = {81}, 
  journal = {Algebra Universalis}, 
  year = {2020}
}

@article{BartoKapytka,
  author    = {Libor Barto and Maryia Kapytka},
  title     = {Multisorted {B}oolean clones determined by binary relations up to minion homomorphisms},
  journal   = {Algebra Universalis},
  volume    = {86(1)},
  year      = {2025},
  url       = {https://doi.org/10.1007/s00012-024-00878-0}}

@Misc{RadekOlsak,
  author = {Radek Ol\v{s}\'{a}k}, 
  title  = {Products of {B}oolean {C}lones up to {M}inion
{H}omomorphisms}, 
  year   = 2024,
  howpublished = {Masters thesis, Charles Univeristy Prague} 
}

@article{BodirskyVucajZhuk,
  AUTHOR = {Bodirsky, Manuel and Vucaj, Albert and Zhuk, Dmitriy},
     TITLE = {The lattice of clones of self-dual operations collapsed},
  JOURNAL = {International Journal of Algebra and Computation},
    VOLUME = {33},
      YEAR = {2023},
    NUMBER = {4},
     PAGES = {717--749},
      ISSN = {0218-1967,1793-6500},
   MRCLASS = {08A40 (08A70)},
  MRNUMBER = {4615620},
MRREVIEWER = {Thodsaporn\ Kumduang},
       DOI = {10.1142/S0218196723500327},
       URL = {https://doi.org/10.1142/S0218196723500327},
}

@book{PosKal79,
 address              = {Berlin},
 author               = {P{\"o}schel, Reinhard and Kalu{\v{z}}nin, Lev A},
 mrclass              = {03G20 (94C10)},
 mrnumber             = {543839 (81f:03075)},
 mrreviewer           = {V. V. Gorlov},
 pages                = {259},
 publisher            = {{VEB} {D}eut\-scher {V}er\-lag der {W}is\-sen\-schaf\-ten},
 series               = {Ma\-the\-ma\-ti\-sche {M}o\-no\-gra\-phi\-en},
 title                = {Funk\-tio\-nen- und {R}e\-la\-tio\-nen\-al\-ge\-bren},
 volume               = {15},
 year                 = 1979,
 URL                  = {https://dx.doi.org/10.1007/978-3-0348-5547-1},
 DOI                  = {10.1007/978-3-0348-5547-1},
 ISBN                 = {978-3-0348-5548-8},
}

@book{Bod21,
    AUTHOR = {Bodirsky, Manuel},
     TITLE = {Complexity of infinite-domain constraint satisfaction},
    SERIES = {Lecture Notes in Logic},
    VOLUME = 52,
 PUBLISHER = {CUP},
   ADDRESS = {Cambridge},
      YEAR = 2021,
     PAGES = {ix+524},
      ISBN = {978-1-107-04284-1},
   MRCLASS = {08-02 (03Cxx 08A70 68Qxx)},
  MRNUMBER = {4273453},
       DOI = {10.1017/9781107337534},
       URL = {https://doi.org/10.1017/9781107337534},
}

@book{McKMcnTay88,
    AUTHOR = {Mc{K}enzie, Ralph N. and
              Mc{N}ulty, George F. and
              Taylor, Walter F.},
     TITLE = {Algebras, lattices, varieties. {V}ol. {I}},
    SERIES = {The Wadsworth \& Brooks/Cole Mathematics Series},
 PUBLISHER = {Wadsworth \& Brooks/Cole Advanced Books \& Software},
   ADDRESS = {Monterey, CA},
    VOLUME = {{I}},
      YEAR = 1987,
     PAGES = {xvi+361},
      ISBN = {0-534-07651-3},
   MRCLASS = {08-01 (06-01)},
  MRNUMBER = {883644 (88e:08001)},
MRREVIEWER = {Gudrun Kalmbach},
}

@book{BS,
  author = "Stanley N. Burris and Hanamantagouda P. Sankappanavar",
  title = "A Course in Universal Algebra",
  publisher = "Springer Verlag, Berlin",
  year = "1981"}

@Article{BoKaKoRo,
  author = {V. G. Bodnar\v{c}uk and L. A. Kalu\v{z}nin and V. N. Kotov and B. A. Romov},
  title = 	 {Galois theory for {Post} algebras, part {I} and {II}},
  journal = 	 {Cybernetics},
  year = 	 1969,
  volume = 5,
  pages = {243--539}
}

@article{Geiger,
title = {Closed Systems of functions and predicates}, 
author = {David Geiger}, 
journal = {Pacific Journal of Mathematics},
volume = {27},
year = {1968},
pages = {95--100}}

@book{FreeseMcKenzie,
author = {Ralph Freese and Ralph McKenzie}, 
title = {Commutator Theory for Congruence Modular Varieties}, 
year = {1987}, 
isbn = {}, 
publisher = {LMS, Lecture Notes (125), Cambridge University Press}
}

@article{AichingerMayrpqextenstions,
  title={Polynomial clones on groups of order pq},
  author={Aichinger, Erhard and Mayr, Peter},
  journal={Acta Mathematica Hungarica},
  volume={114},
  number={3},
  pages={267--285},
  year={2007},
  publisher={Springer}
}

@unpublished{Zeb,
  author = {Brady, Zarathustra},
  title = {Notes on {CSP}s and Polymorphisms},
  publisher = {arXiv},
  year = {2022},
  copyright = {arXiv.org perpetual, non-exclusive license},
  url = {https://arxiv.org/abs/2210.07383},
  note={arXiv preprint arXiv:2210.07383}
}

@article{BakerPixley,
  title={Polynomial interpolation and the Chinese remainder theorem for algebraic systems},
  author={Baker, Kirby A. and Pixley, Alden F.},
  journal={Mathematische Zeitschrift},
  volume={143},
  pages={165--174},
  year={1975},
  publisher={Springer}
}

@book{Bergman,
  author = "Clifford Bergman",
  title = "Universal Algebra: Fundamentals and Selected Topics",
  publisher = "CRC Press",
  year = "2012"}

@article{Bir-On-the-structure,
author={Garrett Birkhoff},
title={On the structure of abstract algebras},
journal={Mathematical Proceedings of the Cambridge Philosophical Society},
volume=31,
number=4,
year=1935,
pages={433--454}
}

@incollection{BagyinskiDemetrovics82,
    AUTHOR = {Demetrovics, J\'{a}nos and Bagyinszki, J\'{a}nos},
     TITLE = {The lattice of linear classes in prime-valued logics},
 BOOKTITLE = {Discrete mathematics},
    SERIES = {{B}anach Center Publications},
    VOLUME = 7,
     PAGES = {105--123},
 PUBLISHER = {PWN--Polish Scientific Publishers},
   ADDRESS = {Warsaw},
      YEAR = 1982,
       DOI = {10.4064/-7-1-105-123},
       URL = {https://doi.org/10.4064/-7-1-105-123},
   MRCLASS = {03G10 (03B50 06B99)},
  MRNUMBER = {698101},
MRREVIEWER = {Daniel Ponasse}
}

@book{HobbyMcKenzie,
author = {David Hobby and Ralph McKenzie},
TITLE = {The structure of finite algebras},
SERIES = {Contemporary Mathematics},
VOLUME = {76},
PUBLISHER = {American Mathematical Society},
YEAR = {1988}
}

@unpublished{MinTaylorOver3,
author = { Libor Barto and Zarathustra Brady and Filip Jankovec and Albert Vucaj and Dmitriy Zhuk},
title = {Minimal {T}aylor clones on three elements},
note = {unpublished},
year = {2025}}

@article{AichingerMayrMcKenzie,
AUTHOR = {Aichinger, Erhard and Mayr, Peter and McKenzie, Ralph},
     TITLE = {On the number of finite algebraic structures},
   JOURNAL = {Journal of the European Mathematical Society (JEMS)},
    VOLUME = {16},
      YEAR = {2014},
    NUMBER = {8},
     PAGES = {1673--1686},
      ISSN = {1435-9855,1435-9863},
   MRCLASS = {08A62 (06A06 08A40 08B05)},
  MRNUMBER = {3262454},
MRREVIEWER = {Neboj\v{s}a\ Mudrinski},
       DOI = {10.4171/JEMS/472},
       URL = {https://doi.org/10.4171/JEMS/472},
}

@unpublished{BM-conservativemaltsev,
  title={{Conservative Maltsev Constraint Satisfaction Problems}},
  author={Bodirsky, Manuel and Moorhead, Andrew},
  note={arXiv preprint arXiv:2505.11395},
  year={2025}
}

@Misc{BulatovMaltsev, 
title = {Malt'sev constraints are tractable},
author = {Andrei A. Bulatov},
howpublished = {Technical report PRG-RR-02-05, Oxford University},
year = {2002}}

\end{document}